\documentclass[10pt]{article}

\usepackage[T1]{fontenc}
\usepackage[english]{babel}
\usepackage{mathrsfs}
\usepackage{amssymb}	
\usepackage{amsthm}
\usepackage{xcolor}
\usepackage{amsmath}
\usepackage{multirow}
\usepackage{nccmath}
\usepackage{rotating}
\usepackage{todonotes}
\usepackage{hyperref}
\usepackage{soul}
\usepackage{empheq}
\usepackage{array,longtable}

\usepackage[left = 2cm, right = 1.5cm, bottom = 5cm, top = 4cm]{geometry}

\newcommand{\R}{\mathbb R}

\newcommand{\weg}[1]{}

\usepackage{graphicx}

\theoremstyle{plain}
\newtheorem{theorem}{Theorem}[section]

\newtheorem{lemma}[theorem]{Lemma}
\newtheorem{proposition}[theorem]{Proposition}

\theoremstyle{definition}

\newtheorem{remark}[theorem]{Remark}

\usepackage{epsf,epsfig,amsfonts,graphicx,caption}

\title{On real-analytic non-degeneracy \\of left-symmetric algebras in dimension 3}
\author{Matteo Manenti\footnote{University of Antwerp, Department of Mathematics, Middelheimlaan 1, B-2020 Antwerpen, Belgium, matteo.manenti@uantwerpen.be}, Serena Scapucci\footnote{Institut f\"ur Mathematik, Friedrich Schiller Universit\"at Jena, 07743 Jena, Germany, serena.scapucci@uni-jena.de}}
\date{}
\usepackage{comment}

\begin{document}
\maketitle

\begin{abstract}
  In this paper, we show that  any three-dimensional differentially non-degenerate left-symmetric algebra over $\mathbb R$ is non-degenerate both in the formal and real-analytic categories via the correspondence with Nijenhuis operators. 
\end{abstract}

\section{Introduction}

Let $\mathfrak{a}$ be an algebra with operation $\star$ over the field $\R$ or $\mathbb C$, then we can naturally define on such an algebra the \textit{associator} $\mathcal{A}$, which is the trilinear operator
\begin{equation*}
	\mathcal{A}(\xi,\eta,\zeta)=(\xi\star\eta)\star\zeta-\xi\star(\eta\star\zeta),\qquad\forall\xi,\eta,\zeta\in\mathfrak{a}.
\end{equation*}
As introduced by Vinberg in \cite{vinberg1967theory} and Koszul in \cite{koszul1961domaines}, an algebra $\mathfrak{a}$ is called \textit{left-symmetric} if $\mathcal{A}(\xi,\eta,\zeta)=\mathcal{A}(\eta,\xi,\zeta)$, (in what follows, we will write LSA for left-symmetric algebra). Though originally introduced in the study of homogeneous domains and convex cones, these algebras have since appeared in many other contexts. In particular, one can tackle the study of such algebras from the point of view of Nijenhuis geometry, thanks to the following construction.\\

Fix any $n$-dimensional left-symmetric algebra $\mathfrak{a}$ and an element $\eta\in\mathfrak{a}$. Then, the algebra admits a natural affine space structure, i.e.\ introducing a basis $\eta_1,\dots,\eta_n$ we have natural coordinates $\eta=x_1 \eta_1+\dots+x_n \eta_n$. We can then consider the operator field, i.e. a $(1,1)$-tensor field, given by
\begin{equation*}
    (R_\mathfrak{a})_\eta(\xi)=\xi\star\eta,\qquad \xi\in T_\eta\mathfrak{a}\cong\mathfrak{a}.
\end{equation*}
For simplicity, we omit $\eta$ and denote it by $R_\mathfrak{a}$. This operator is called \textit{the right-adjoint} operator of $\mathfrak{a}$, and in coordinates $x_1,\dots,x_n$, the corresponding matrix has entries
$(R_\mathfrak{a})^i_j=\sum_{k=1}^n a^i_{jk}x_k$,
where $a^i_{jk}$ are the structure constants of the algebra. Thus, all the entries of the right-adjoint operator are linear in the coordinates. 

\begin{theorem}[\cite{winterhalder1997linear}]
	Let $(\mathfrak{a},\star)$ be an algebra over $\R$ or $\mathbb C$ and let $R_{\mathfrak{a}}$ be the right-adjoint operator field. Then $\mathfrak{a}$ is a left symmetric algebra if and only if the Nijenhuis torsion of $R_\mathfrak{a}$
	\begin{equation}\label{eq: Nij torsion}
		\mathcal{N}_{R_\mathfrak{a}}(X,Y):=R_\mathfrak{a}^2[X,Y]+[R_\mathfrak{a}X,R_\mathfrak{a}Y]-R_\mathfrak{a}[R_\mathfrak{a}X,Y]-R_\mathfrak{a}[X,R_\mathfrak{a}Y]
	\end{equation}
	vanishes for all vector fields $X,Y$.
\end{theorem}

The Nijenhuis torsion of an operator field $R$ on a manifold $M^n$ was first defined in \cite{nijenhuis1951xn}, and if $\mathcal{N}_R\equiv0$ we say that the operator $R$ is \textit{Nijenhuis}. As already pointed out in \cite{konyaev2024linearizationcertainsingularitiesnijenhuis}, the above result establishes the analogous role of left-symmetric algebras in Nijenhuis geometry as the one of Lie algebras in Poisson geometry. In the framework of Poisson geometry, the problem of linearization of Poisson brackets vanishing at a point was introduced by Weinstein in \cite{weinstein1983local}, and following the same path in \cite{bolsinov2022nijenhuisNGI} the problem of linearization of Nijenhuis operators was proposed.\\

Consider a manifold $M^n$ and a Nijenhuis operator $R$ on it. Assume that $R(p)=\lambda_0\cdot\operatorname{Id}$ for some $p\in M^n$ and for some constant $\lambda_0$. In this case, we call $p$ a \textit{singular point of scalar type}. Fixing local coordinates $x = (x_1,\dots,x_n)$ centered at $p$ we Taylor expand $R$ and get $R(x)=\lambda_0\cdot\operatorname{Id}+R_1(x)+R_2(x)+\dots$, where the entries of $R_1(x)$ are linear in $x_1,\dots,x_n$. Then, $R_1$ is itself Nijenhuis and it induces a natural structure of left-symmetric algebra on $T_pM^n$. 
We say that $R$ is \textit{linearizable} (in the formal, resp. real-analytic category) if there exists a coordinate change such that $R=\lambda_0\cdot\operatorname{Id}+R_1$. \\

On the algebras side, we say that a left-symmetric algebra $\mathfrak{a}$ is \textit{non-degenerate} in the formal (resp. real-analytic or smooth) category if any Nijenhuis operator with linear part $R_1=R_{\mathfrak{a}}$ is linearizable in the formal (resp. real-analytic or smooth) category. While in dimension 2 the study of non-degeneracy of left-symmetric algebras has been completed in \cite{konyaev2021nijenhuisNGII,konyaev2024linearizationcertainsingularitiesnijenhuis}, in higher dimension the problem is still open. In fact, there is not even a full classification of such algebras in dimension $\geq3$ over the reals, even though a classification of complex three-dimensional LSAs has been found in \cite{Bai_2009}. Nevertheless, in dimension 3 we have some classification results of families of Nijenhuis operators satisfying stronger properties in the following sense. \\

A point $p\in M$ is called \textit{differentially non-degenerate} for an operator $R$, if the differentials d$\sigma_1(p),\dots,$d$\sigma_n(p)$ of the coefficients of the characteristic polynomial $\chi_{R(p)}(t)=\sum_{k=0}^{n}t^k\sigma_{n-k}(x)$ are linearly independent at this point. In this case we say that the functions $\sigma_1,\dots,\sigma_n$ are \textit{functionally independent} at the point $p\in M$, and that the operator $R$ is \textit{differentially non-degenerate} at the point $p\in M$. We say that a LSA $\mathfrak{a}$ is \textit{differentially non-degenerate} if the corresponding right-adjoint Nijenhuis operator $R_\mathfrak{a}$ is differentially non-degenerate almost everywhere.\\

Differentially non-degenerate LSAs in dimension 3 have been fully classified in \cite{scapucci2026classification,degtyareva2024classification}, and we recall this classification in the following Theorem.

\begin{theorem}\label{thm:LSA3d}
	The classification in dimension 3 of differentially non-degenerate LSAs over $\R$ or equivalently Nijenhuis operators with linear entries whose coefficients of the characteristic polynomial are functionally independent at almost every point is presented in Table 1 and Table 2, the decomposable and indecomposable ones respectively\footnote{We say that a LSA $\mathfrak{a}$ is \textit{decomposable} over a field $\mathbb K$ if it can be written as a direct product of lower dimensional LSAs. We say that the LSA is \textit{indecomposable} otherwise.}\footnote{The notations for the decomposable algebras are taken from \cite{bolsinov2022nijenhuisNGI} and \cite{konyaev2021nijenhuisNGII}: $\mathfrak{d}$ denotes the non-trivial one-dimensional LSA, the $\mathfrak{b}$-series represent the two-dimensional LSAs with non-commutative associated Lie algebra and finally the $\mathfrak{c}$-series represent the two-dimensional LSAs with commutative associated Lie algebra.}. For each case the structure relations of the LSA, the corresponding right-adjoint Nijenhuis operator and the coefficients of its characteristic polynomial are given.
	
	\begingroup
	\small
	\setlength{\tabcolsep}{6pt}
	\setlength{\LTleft}{\fill}   
	\setlength{\LTright}{\fill} 
	
	\begin{longtable}{
			|>{\centering\arraybackslash}m{0.08\textwidth}
			|>{\centering\arraybackslash}m{0.22\textwidth}
			|>{\centering\arraybackslash}m{0.32\textwidth}
			|>{\centering\arraybackslash}m{0.28\textwidth}|
		}
		\caption{Decomposable 3-dimensional LSAs}
		\label{tab:decomposable-lsa}\\
		\hline
		\rule[-1.8ex]{0pt}{4.5ex}\textbf{Name}
		&
		\rule[-1.8ex]{0pt}{4.5ex}\textbf{Structure relations}
		&
		\rule[-1.8ex]{0pt}{4.5ex}\textbf{Nijenhuis operator $R_{\mathfrak{a}}$}
		&
		\rule[-1.8ex]{0pt}{4.5ex}\textbf{Coefficients of $\chi_{R_{\mathfrak{a}}}$}
		\tabularnewline
		\hline
		\endfirsthead
		
		\multicolumn{4}{c}{\tablename\ \thetable\ (continued)}\\
		\hline
		\textbf{Name}
		&
		\textbf{Structure relations}
		&
		\textbf{Nijenhuis operator $R_{\mathfrak{a}}$}
		&
		\textbf{Coefficients of $\chi_{R_{\mathfrak{a}}}$}
		\tabularnewline
		\hline
		\endhead
		
		\hline
		\multicolumn{4}{r}{\emph{Continued on the next page}}\\
		\endfoot
		
		\hline
		\endlastfoot
		
		\rule[-3ex]{0pt}{7.5ex} $\mathfrak{b_4^+}\oplus\mathfrak{d}$
		&
		$\begin{aligned}
			\eta_1\star\eta_1 &= 2\eta_1 \\
			\eta_1\star\eta_2 &= \eta_2 \\
			\eta_2\star\eta_2 &= -\eta_1 \\
			\eta_3\star\eta_3 &= \eta_3
		\end{aligned}$
		&
		$\left(\begin{array}{ccc}
			2x_1 & -x_2 & 0\\
			x_2 & 0 & 0\\
			0 & 0 & x_3
		\end{array}\right)$
		&
		$\begin{aligned}
			\sigma_1 &= -2x_1-x_3\\
			\sigma_2 &= x_2^2+2x_1x_3\\
			\sigma_3 &= -x_2^2x_3
		\end{aligned}$
		\tabularnewline
		\hline
		
		\rule[-3ex]{0pt}{7.5ex} $\mathfrak{b_4^-}\oplus\mathfrak{d}$
		&
		$\begin{aligned}
			\eta_1\star\eta_1 &= 2\eta_1 \\
			\eta_1\star\eta_2 &= \eta_2 \\
			\eta_2\star\eta_2 &= \eta_1 \\
			\eta_3\star\eta_3 &= \eta_3
		\end{aligned}$
		&
		$\left(\begin{array}{ccc}
			2x_1 & x_2 & 0 \\
			x_2 & 0 & 0 \\
			0 & 0 & x_3
		\end{array}\right)$
		&
		$\begin{aligned}
			\sigma_1 &= -2x_1-x_3\\
			\sigma_2 &= -x_2^2+2x_1x_3\\
			\sigma_3 &= x_2^2x_3
		\end{aligned}$
		\tabularnewline
		\hline
		
		\rule[-3ex]{0pt}{7.5ex} $\mathfrak{c_5^+}\oplus\mathfrak{d}$
		&
		$\begin{aligned}
			\eta_1\star\eta_1 &= \eta_1 \\
			\eta_1\star\eta_2 &= \eta_2 \\
			\eta_2\star\eta_1 &= \eta_2 \\
			\eta_2\star\eta_2 &= \eta_1 \\
			\eta_3\star\eta_3 &= \eta_3
		\end{aligned}$
		&
		$\left(\begin{array}{ccc}
			x_1 & x_2 & 0 \\
			x_2 & x_1 & 0 \\
			0 & 0 & x_3
		\end{array}\right)$
		&
		$\begin{aligned}
			\sigma_1 &= -2x_1-x_3\\
			\sigma_2 &= x_1^2-x_2^2+2x_1x_3\\
			\sigma_3 &= -x_1^2x_3+x_2^2x_3
		\end{aligned}$
		\tabularnewline
		\hline
		
		\rule[-3ex]{0pt}{7.5ex} $\mathfrak{c_5^-}\oplus\mathfrak{d}$
		&
		$\begin{aligned}
			\eta_1\star\eta_1 &= \eta_1 \\
			\eta_1\star\eta_2 &= \eta_2 \\
			\eta_2\star\eta_1 &= \eta_2 \\
			\eta_2\star\eta_2 &= -\eta_1 \\
			\eta_3\star\eta_3 &= \eta_3
		\end{aligned}$
		&
		$\left(\begin{array}{ccc}
			x_1 & -x_2 & 0 \\
			x_2 & x_1 & 0 \\
			0 & 0 & x_3
		\end{array}\right)$
		&
		$\begin{aligned}
			\sigma_1 &= -2x_1-x_3\\
			\sigma_2 &= x_1^2+x_2^2+2x_1x_3\\
			\sigma_3 &= -x_1^2x_3-x_2^2x_3
		\end{aligned}$
		\tabularnewline
		\hline
	\end{longtable}
	\endgroup
	
	\begingroup
	\small
	\setlength{\tabcolsep}{6pt}
	\setlength{\LTleft}{\fill}  
	\setlength{\LTright}{\fill}
	
	\begin{longtable}{
			|>{\centering\arraybackslash}m{0.08\textwidth}
			|>{\centering\arraybackslash}m{0.22\textwidth}
			|>{\centering\arraybackslash}m{0.32\textwidth}
			|>{\centering\arraybackslash}m{0.28\textwidth}|
		}
		\caption{Indecomposable 3-dimensional LSAs}
		\label{tab:indecomposable-lsa}\\
		\hline
		\rule[-1.8ex]{0pt}{4.5ex}\textbf{Name}
		&
		\rule[-1.8ex]{0pt}{4.5ex}\textbf{Structure relations}
		&
		\rule[-1.8ex]{0pt}{4.5ex}\textbf{Nijenhuis operator $R_{\mathfrak{a}}$}
		&
		\rule[-1.8ex]{0pt}{4.5ex}\textbf{Coefficients of $\chi_{R_{\mathfrak{a}}}$}
		\tabularnewline
		\hline
		\endfirsthead
		
		\multicolumn{4}{c}{\tablename\ \thetable\ (continued)}\\
		\hline
		\textbf{Name}
		&
		\textbf{Structure relations}
		&
		\textbf{Nijenhuis operator $R_{\mathfrak{a}}$}
		&
		\textbf{Coefficients of $\chi_{R_{\mathfrak{a}}}$}
		\tabularnewline
		\hline
		\endhead
		
		\hline
		\multicolumn{4}{r}{\emph{Continued on the next page}}\\
		\endfoot
		
		\hline
		\endlastfoot
		
		\rule[-3ex]{0pt}{7.5ex} $\mathfrak{s}_1$
		&
		$\begin{aligned}
			\eta_1\star\eta_1 &= 2\eta_1 \\
			\eta_1\star\eta_2 &= \eta_2 \\
			\eta_1\star\eta_3 &= -\eta_1 \\
			\eta_2\star\eta_2 &= -\eta_1 \\
			\eta_2\star\eta_3 &= \eta_2 \\
			\eta_3\star\eta_1 &= -\eta_1 \\
			\eta_3\star\eta_3 &= \eta_1 + \eta_3
		\end{aligned}$
		&
		$\left(\begin{array}{ccc}
			2x_1-x_3 & -x_2 & x_3-x_1 \\
			x_2 & x_3 & 0 \\
			0 & 0 & x_3
		\end{array}\right)$
		&
		$\begin{aligned}
			\sigma_1 &= -2x_1-x_3\\
			\sigma_2 &= x_2^2 + 4x_1x_3 - x_3^2\\
			\sigma_3 &= - x_2^2x_3 - 2x_1x_3^2 + x_3^3
		\end{aligned}$
		\tabularnewline
		\hline
		
		\rule[-3ex]{0pt}{7.5ex} $\mathfrak{s}_2$
		&
		$\begin{aligned}
			\eta_1\star\eta_1 &= 2\eta_1 \\
			\eta_1\star\eta_2 &= \eta_2 \\
			\eta_1\star\eta_3 &= -\eta_1 \\
			\eta_2\star\eta_2 &= \eta_1 \\
			\eta_2\star\eta_3 &= \eta_2 \\
			\eta_3\star\eta_1 &= -\eta_1 \\
			\eta_3\star\eta_3 &= \eta_1 + \eta_3
		\end{aligned}$
		&
		$\left(\begin{array}{ccc}
			2x_1-x_3 & x_2 & x_3-x_1 \\
			x_2 & x_3 & 0 \\
			0 & 0 & x_3
		\end{array}\right)$
		&
		$\begin{aligned}
			\sigma_1 &= -2x_1-x_3 \\
			\sigma_2 &= -x_2^2 + 4x_1x_3 - x_3^2\\
			\sigma_3 &= x_2^2x_3 - 2x_1x_3^2 + x_3^3
		\end{aligned}$
		\tabularnewline
		\hline
		
		\rule[-3ex]{0pt}{7.5ex} $\mathfrak{s}_3$
		&
		$\begin{aligned}
			\eta_1\star\eta_1 &= \eta_1 \\
			\eta_1\star\eta_2 &= \eta_2 \\
			\eta_2\star\eta_3 &= -\eta_1 \\
			\eta_3\star\eta_2 &= -\eta_1 -\eta_2\\
			\eta_3\star\eta_3 &= \eta_3
		\end{aligned}$
		&
		$\left(\begin{array}{ccc}
			x_1 & -x_3 & -x_2 \\
			x_2 & 0 & -x_2 \\
			0 & 0 & x_3
		\end{array}\right)$
		&
		$\begin{aligned}
			\sigma_1 &= -x_1-x_3 \\
			\sigma_2 &= x_1x_3 + x_2x_3 \\
			\sigma_3 &= -x_2x_3^2
		\end{aligned}$
		\tabularnewline
		\hline
		
		\rule[-3ex]{0pt}{7.5ex} $\mathfrak{s}_4$
		&
		$\begin{aligned}
			\eta_1\star\eta_1 &= -\eta_1 \\
			\eta_1\star\eta_2 &= -\tfrac{2}{3}\eta_2 \\
			\eta_1\star\eta_3 &= -\tfrac{1}{3}\eta_3 \\
			\eta_2\star\eta_3 &= \eta_1 \\
			\eta_3\star\eta_2 &= \eta_1 \\
			\eta_3\star\eta_3 &= \eta_2
		\end{aligned}$
		&
		$\left(\begin{array}{rrr}
			-x_1 & x_3 & x_2 \\
			-\frac{2}{3} \, x_2 & 0 & x_3 \\
			-\frac{1}{3} \, x_3 & 0 & 0
		\end{array}\right)$
		&
		$\begin{aligned}
			\sigma_1 &= x_1 \\
			\sigma_2 &= x_2x_3 \\
			\sigma_3 &= \frac{1}{3}x_3^3
		\end{aligned}$
		\tabularnewline
		\hline
	\end{longtable}
	\endgroup  
\end{theorem}

Therefore, it is a natural question whether the algebras in this classification are non-degenerate. The present paper is devoted to study this problem. Specifically, we establish the first non-degeneracy results for three-dimensional left-symmetric algebras beyond the diagonal case, previously settled in \cite[Theorem 5.1]{bolsinov2022nijenhuisNGI}. Namely, the main result of this paper is the following.

\begin{theorem}\label{thm: main result introduction}
	Any 3-dimensional differentially non-degenerate left-symmetric algebra is non-degenerate in both the formal and real-analytic categories.
\end{theorem}

This result fits into the study of singular points of Nijenhuis operators, which is an important part of the Nijenhuis geometry programme, initiated in \cite{bolsinov2022nijenhuisNGI}. As one can notice already in this seminal paper, differentially non-degenerate Nijenhuis operators arise naturally in many applications of Nijenhuis geometry, and in particular they play an important role when dealing with problems related to integrable systems. \\

For instance, in \cite{bolsinov2022applications}, the authors show how differentially non-degenerate operators allow one to construct series of multi-component integrable PDEs with differential constraints, which generalize many famous integrable systems including KdV, Camassa-Holm and others (see \cite[Theorem 2]{bolsinov2022applications}). Furthermore, differentially non-degenerate operators appear naturally in the study of Dubrovin-Frobenius geometry. Indeed, a differentially non-degenerate Nijenhuis operator $L$ always admits a \textit{unity}, i.e.\ a vector field $e$ such that $\mathcal L_eL=\operatorname{Id}$, and the class of Nijenhuis manifolds admitting a cyclic unity coincides with regular $F$-manifolds (see \cite{antonov2024nijenhuis}). The unity of a differentially non-degenerate operator is cyclic outside the singular set of the operator. Thus, any differentially non-degenerate operator is associated naturally with an F-manifold structure on the regular set of the operator.\\

Moreover, it is known that differential non-degeneracy of a Nijenhuis operator at a point implies $\operatorname{gl}$-regularity at  the same point (i.e.\ its Jordan normal form contains exactly one Jordan block for each eigenvalue). In particular, differentially non-degenerate operators form a subclass of $\operatorname{gl}$-regular operators, which are the natural framework for constructing integrable systems of hydrodynamic type, as shown in \cite{bolsinov2024nijenhuis,bolsinov2024applications}. They have also been applied to the theory of (generalized) separation of variables (\cite{bolsinov2022applications,bolsinov2025orthogonal,bolsinov2025st}). These applications further motivate the study of the local structure of Nijenhuis operators, and in particular of their behaviour near singular points, especially of scalar type, where the operator is as far as possible from being $\operatorname{gl}$-regular.\\

The main result of this paper also has potential applications in the classical theory of geodesically equivalent metrics and projective vector fields, in view of \cite[Theorem 6.2]{bolsinov2022nijenhuisNGI}, where the authors exploit differential non-degeneracy of an operator $L$ in a point $p$ to describe \textit{all} the metrics compatible with the operator in the real-analytic category. Indeed, Theorem \ref{thm: main result introduction} represents a first tool to partially tackle the Beltrami problem in dimension 3. Recently, Akpan and Bolsinov solved the Beltrami problem of describing geodesically equivalent metrics in dimension 2 by exploiting the description of singular points of the corresponding Nijenhuis operator (see \cite{akpan2026beltrami}), and hence we expect that our result can be used in a similar way near points of scalar type for three-dimensional metrics.\\

Typically, all these problems rely on the fact that the description of the corresponding geometric/analytic objects in neighbourhoods of differentially non-degenerate points is well known, and it depends on certain parameters. To study these problems, one needs to understand (case by case) how to choose the parameters in order to extend the structures to singular points. In this sense, the classification of \cite{scapucci2026classification} and the linearization (Theorem \ref{thm: main result introduction}) provide a starting point to deal with similar questions in dimension three.\\

It is worth noticing that Theorem \ref{thm: main result introduction} describes linearization in the formal and real-analytic categories, but it does not give any information in the smooth category. Furthermore, in general the following problem remains open: given two non-degenerate LSAs, is their direct product still non-degenerate? This problem was already proposed in \cite{bkm2025}, and the present result shows that this is true when we consider differentially non-degenerate LSAs in dimension three. We plan to tackle these further problems and some of the questions described above in future works.

\section*{Acknowledgements}
The authors express their gratitude to A.\ V.\ Bolsinov, V.\ S.\ Matveev and D.\ Akpan for their attention to the text and for valuable discussions. Part of the project was carried out during the ``Poisson Summer School 2026", held at the University of Antwerp, and it was inspired by the problems presented in \cite{Bolsinov2026}.

\section*{Funding}
The first author is supported by the FWO-EoS grant through the project \textit{Beyond symplectic geometry} with AU 45816. The second author is partially funded by DFG, project 529233771. 

\newpage
\section{Main Result}
Fixed the notation as above, we reformulate Theorem \ref{thm: main result introduction} in terms of Nijenhuis operators.

\begin{theorem}\label{thm:main result}
    Let $\mathfrak{a}$ be any of the left-symmetric algebras of Theorem \ref{thm:LSA3d}. Then, $\mathfrak{a}$ is non-degenerate in both the formal and real-analytic categories.\\
    Equivalently, consider the right-adjoint Nijenhuis operator $R_{\mathfrak{a}}$ corresponding to $\mathfrak{a}$ and a formal (resp., real analytic) perturbation of $R_{\mathfrak{a}}$ of the form
    \begin{equation}\label{eq: perturbation0}
    R(x)=R_{\mathfrak{a}}(x)+R_2(x)+R_3(x)+\dots,
    \end{equation}
    where the entries of $R_k(x)$ are homogeneous polynomials of degree $k \ge 2$ in $x$ and $R$ is a formal (resp. real analytic) Nijenhuis operator. Then $R(x)$ is formally (resp., real analytically) linearizable, i.e., there exists a formal (resp. real analytic) change of variables 
    \begin{equation}\label{eq: coordinate change}
        y_j=x_j+\psi_{j,2}(x)+\psi_{j,3}(x)+\dots
    \end{equation}
    where $\psi_{j,k}$ are homogeneous polynomials of degree $k \ge 2$ such that in new coordinates the entries of the operator $R$ are linear, or, more specifically, 
    \begin{equation}\label{eq: formal lin}
    \left(\frac{\partial y}{\partial x}\right)\left(R_{\mathfrak{a}}(x)+R_2(x)+R_3(x)+\dots\right)\left(\frac{\partial y}{\partial x}\right)^{-1} =R_{\mathfrak{a}}(y),
    \end{equation}
    where $\left(\frac{\partial y}{\partial x}\right)$ is the Jacobi matrix of the coordinate change \eqref{eq: coordinate change}.
\end{theorem}

\subsection{Proof}\label{sec: approach}
We fix a left-symmetric algebra $\mathfrak{a}$ with associated right-adjoint Nijenhuis operator $R_{\mathfrak{a}}$. We prove Theorem \ref{thm:main result} showing first formal linearizability and then we apply Artin's Theorem to get real-analytic linearizability.

\subsubsection*{Formal linearizability}

Consider a Nijenhuis formal perturbation of $R_{\mathfrak{a}}$ of the form 
\begin{equation}\label{eq: perturbation}
	R(x)=R_{\mathfrak{a}}(x)+R_k(x)+R_{k+1}(x)+\dots,
\end{equation}
where
\begin{equation}\label{eq: k-component_of_pert}
R_k=\begin{pmatrix}
        f_{11}&f_{12}&f_{13}\\
        f_{21}&f_{22}&f_{23}\\
        f_{31}&f_{32}&f_{33}
    \end{pmatrix}
\end{equation}
with entries $f_{ij}$ that are homogeneous polynomials of degree $k$. Notice that if $R_\mathfrak{a}$ is differentially non-degenerate, then $R$ as in \eqref{eq: perturbation} is still differentially non-degenerate.\\

Recall that given any operator fields $M,N$, their \textit{Fr\"olicher-Nijenhuis brackets} are
\begin{equation*}
    [M,N]_{FN}:=\mathcal{N}_{M+N}-\mathcal{N}_M-\mathcal{N}_N
\end{equation*}
Now, we notice that if a perturbation of the form \eqref{eq: perturbation} is Nijenhuis, then the Fr\"olicher-Nijenhuis bracket $[R_{\mathfrak{a}},R_{k}]_{FN}$ vanishes. Indeed, since $R$ is Nijenhuis, $[R,R]_{FN}$ has to vanish identically and $[R_{\mathfrak{a}},R_{k}]_{FN}$ is the homogeneous part of degree $k$ of
\begin{align*}
    [R,R]_{FN}&=[R_{\mathfrak{a}}+R_k+R_{k+1}+\dots,R_{\mathfrak{a}}+R_k+R_{k+1}+\dots]_{FN}\\&=[R_{\mathfrak{a}},R_{\mathfrak{a}}]_{FN}+2[R_{\mathfrak{a}},R_k]_{FN}+[R_k,R_k]_{FN}+2[R_{\mathfrak{a}},R_{k+1}]_{FN}+\dots
\end{align*} 
We then prove formal linearizability by killing each term $R_k(x)$ of the perturbation step by step, making use of the following result.
\begin{proposition}[{\cite[Proposition 5.3]{bolsinov2022nijenhuisNGI}}]\label{prop:BKM1,5.3}
    Assume that for any $k$ and for any $R_k$ such that $[R_{\mathfrak{a}},R_k]_{FN}=0$ there
exists a change of variables
$y_i = x_i + \psi_{i,k}(x)$, where $\psi_{i,k}(x)$ is a homogeneous polynomial in $x$ of degree $k$,
that transforms $ R_{\mathfrak{a}}(x) + R_k (x)$ to the form $ R_{\mathfrak{a}}(y) + R_{k+1}(y) +\dots$(i.e.,
kills terms of order $k$), that is
\begin{equation}\label{eq:linearization}
        (\operatorname{Id}+J)^{-1}\left(R_{\mathfrak{a}}(y)+(degree\,k+1)\right)(\operatorname{Id}+J)=R_{\mathfrak{a}}(x)+R_k(x),
\end{equation} 
where $J$ is the Jacobian of the map $\psi_k:=(\psi_{1,k},\psi_{2,k},\psi_{3,k})$. Then any Nijenhuis perturbation of $R_{\mathfrak{a}}$ of the form \eqref{eq: perturbation}
is formally linearizable.
\end{proposition}
To simplify the notation, we write $\psi_j$ instead of $\psi_{j,k}$, since $k$ is fixed in what follows. With the above proposition, proving formal linearizability reduces to finding a coordinate change of the form $y_i = x_i + \psi_{i}(x)$ such that \eqref{eq:linearization} holds. By equating the $k$-order terms of \eqref{eq:linearization} we get
\begin{equation}\label{eq: condition lin}
    [R_{\mathfrak{a}},J]=R_k-R_{\mathfrak{a}}(\psi(x)).
\end{equation}
This system is easy to analyze, and solving it we get the unique coordinate change needed for formal linearizability and conditions on the entries of the operator $R_k$.

\begin{lemma}\label{lemma: main}
    Let $\mathfrak{a}$ be one of the left-symmetric algebras of Theorem \ref{thm:LSA3d}.
    \begin{enumerate}
    \item The change of variables $y_i=x_i+\psi_i(x)$ with  required conditions \eqref{eq: condition lin} exists if and only if $R_k$ satisfies suitable compatibility conditions, which we denote by $(\star)$. Moreover, whenever such a solution exists, it is unique.
    \item The vanishing of the Fr\"olicher-Nijenhuis brackets $[R_{\mathfrak{a}},R_k]_{FN}$ implies
that the relations $(\star)$ are satisfied. In particular, the coordinate change $y_i=x_i+\psi_i(x)$ is well defined and satisfies~\eqref{eq: condition lin}.
    \end{enumerate}
\end{lemma}

\begin{remark}
  The conditions $(\star)$ depend on the algebra $\mathfrak a$. For each algebra, these conditions are given explicitly in the first lemma of the corresponding subsection of Section~\ref{sec: main Lemma}, which will be devoted to the proof of Lemma \ref{lemma: main}. We will study it case by case, describing explicitly the relations $(\star)$ and concluding \textit{formal} linearizability of the operator $R_{\mathfrak{a}}$. 
\end{remark}

\subsubsection*{Application of Artin's Theorem}

Now, in order to show real-analytic linearizability from the formal one, we make use of the following result. 

\begin{theorem}[Artin, \cite{artin1968solutions}]\label{thm: artin}
Let $F_1,\ldots,F_m$ be polynomials in variables
$y_1,\ldots,y_N$, with coefficients that are convergent power
series in $x=(x_1,\ldots,x_n)$. If there exist formal power series
$\hat y_1(x),\ldots,\hat y_N(x)$ such that
\[
\begin{cases}
F_1(x,\hat y(x))=0,\\
\vdots\\
F_m(x,\hat y(x))=0,
\end{cases}
\]
where $\hat{y}(x)=(\hat y_1(x),\ldots,\hat y_N(x))$, then for any $c$ there exist convergent power series
$y_1(x),\ldots,y_N(x)$ such that
\[
\begin{cases}
F_1(x,y(x))=0,\\
\vdots\\
F_m(x,y(x))=0,
\end{cases}
\qquad\text{and}\qquad
y_j(x)\equiv \hat y_j(x)\pmod{(x)^c}.
\]
\end{theorem}

Now we assume that the perturbation $R(x)$ \eqref{eq: perturbation} is real analytic. Let the characteristic polynomial of $R$ be 
\begin{equation*}
    \chi_R(t)=t^3+\zeta_1(x)t^2+\zeta_2(x)t+\zeta_3(x)
\end{equation*}
where the $\zeta_is$ are real analytic functions. \\

By the construction above, there exist formal coordinates $y_j(x)$ such that \eqref{eq: formal lin} holds. In particular, the characteristic polynomials of the two sides coincide, that is, denoting by $\sigma_j$ the coefficients of the characteristic polynomial of $R_{\mathfrak{a}}$, 

\begin{equation*}
    \begin{cases}
        \sigma_1(y)=\zeta_1(x)\\
        \sigma_2(y)=\zeta_2(x)\\
        \sigma_3(y)=\zeta_3(x)
    \end{cases}
\end{equation*}
Notice that $\sigma_j$ are polynomials in $y_1,y_2,y_3$, and $\zeta_j$ are real-analytic functions in $x$. Thus, this system satisfies the hypothesis of Artin's Theorem, and it admits a formal solution, thanks to the paragraph above. Hence, Theorem \ref{thm: artin} guarantees the existence of a real-analytic solution to the system. Now, thanks to \cite[Corollary 2.2]{bolsinov2022nijenhuisNGI}, Nijenhuis operators with the same characteristic polynomials and corresponding (almost everywhere) functionally independent coefficients coincide, then the real-analytic solution of the previous system is a solution also of the equation \eqref{eq: formal lin}, and hence we get the real-analytic linearizability we claimed. 

\begin{remark}
    Notice that this argument is independent on the dimension: any right-adjoint operator for a LSA has linear entries, and the coefficients of the characteristic polynomials are (homogeneous) polynomials of the coordinates, thus the previous argument can be easily generalized to any dimension. 
\end{remark}

\begin{proposition}\label{prop: formal_to_real_an}
    Let $\mathfrak{a}$ be any differentially non-degenerate LSA. If $\mathfrak{a}$ is formally non-degenerate, then $\mathfrak{a}$ is real-analytically non-degenerate.
\end{proposition}

Thus, in what follows, we consider each left symmetric algebra from the list separately, and for each of them, it is enough to prove formal linearizability, because by Proposition \ref{prop: formal_to_real_an} we get automatically the real-analytic linearizability. 

\section{Proof of Lemma \ref{lemma: main}}\label{sec: main Lemma}

In this section, we will prove Lemma \ref{lemma: main} by studying each case separately. In particular, let $R_\mathfrak{a}$ be one of the operators from Theorem \ref{thm:LSA3d} and let $R_k=(f_{ij})$ be an operator such that $[R_{\mathfrak{a}},R_k]_{FN}=0$, with $f_{ij}$ homogeneous polynomials of degree $k$. For each algebra, we will first describe the relations $(\star)$ of the first point of Lemma \ref{lemma: main} and the corresponding coordinate change, by expanding \eqref{eq: condition lin}. This is the content of the first lemma of each subsection. The second Lemma is devoted to the proof of the second point of Lemma \ref{lemma: main}, i.e.\ we prove that $[R_{\mathfrak{a}},R_k]_{FN}=0$ implies $(\star)$  and the coordinate change is well defined.

\subsection{Non-degeneracy of the LSA $\mathfrak{c}_5^{\pm}\oplus\mathfrak{d}$}

As noticed in \cite[Remark 3.1]{scapucci2026classification}, it is known that the LSA $\mathfrak{c}_5^+\oplus\mathfrak{d}\cong\oplus_3\mathfrak{d}$ over $\R$, thus:

\begin{proposition}[{\cite[Theorem 5.1]{bolsinov2022nijenhuisNGI}}]
    Let $\mathfrak{c}_5^+$ be the 2-dimensional LSA associated with linear operator
    \begin{equation*}
        R_{\mathfrak{c}_5^+}=\begin{pmatrix}
            x_1&x_2\\x_2&x_1
        \end{pmatrix}
    \end{equation*}
    and let $\mathfrak{d}$ be a non-trivial one-dimensional LSA. Then $\mathfrak{c}_5^+\oplus\mathfrak{d}$ is non-degenerate in both the formal and real-analytic categories.
\end{proposition}

Let's now consider instead the LSA $\mathfrak{c}_5^-\oplus\mathfrak{d}$, which is associated with the Nijenhuis operator 
\begin{equation}\label{eq: op c5-d}
    R_{\mathfrak{c}_5^-\oplus\mathfrak{d}}=\begin{pmatrix}
        x_1&-x_2&0\\
        x_2&x_1&0\\
        0&0&x_3
    \end{pmatrix}
\end{equation}

\begin{lemma}
    The change of coordinates $y_i=x_i+\psi_i(x)$ with required condition \eqref{eq: condition lin} exists if and only if the components of $R_k$ satisfy the following relations:
    \begin{subequations}\label{eq: conditions for existence c5-}
        \begin{empheq}[left=\empheqlbrace]{align}
            &f_{22}-f_{11}=x_2\left(\partial_1f_{21}-\partial_1f_{12}+\partial_2f_{11}+\partial_2f_{22}\right)\label{eq: c5-,1}\\
            &f_{21}+f_{12}=x_2\left(\partial_1f_{11}+\partial_1f_{22}-\partial_2f_{21}+\partial_2f_{12}\right)\label{eq: c5-,2}\\
            &2 f_{13}=(x_1-x_3)(\partial_3 f_{11}+\partial_3 f_{22})+x_2(\partial_3f_{12}-\partial_3f_{21})\label{eq: c5-,3}\\
            &2f_{23}= (x_1-x_3)(\partial_3f_{21}-\partial_3 f_{12})+x_2(\partial_3f_{11}+\partial_3f_{22})\label{eq: c5-,4}\\
            &f_{31}=(x_3-x_1)\partial_1f_{33}-x_2\partial_2f_{33}\label{eq: c5-,5}\\
            &f_{32}=(x_3-x_1)\partial_2 f_{33}+x_2\partial_1f_{33}\label{eq: c5-,6}
        \end{empheq}
    \end{subequations}
    In this case, the coordinate change takes the form
    \begin{equation*}
        \begin{cases}\vspace{0.1cm}
            y_1=x_1+\dfrac{f_{11}+f_{22}}{2}\\\vspace{0.1cm}
            y_2=x_2+\dfrac{f_{21}-f_{12}}{2}\\
            y_3=x_3+f_{33}
        \end{cases}
    \end{equation*}
\end{lemma}

\begin{proof}
    For the operator \eqref{eq: op c5-d}, equation \eqref{eq: condition lin} is equivalent to
    \begin{subequations}\label{eq: conditions proof c5-}
        \begin{empheq}[left=\empheqlbrace]{align}
            &f_{11}-\psi_1=-x_2(\partial_1\psi_2+\partial_2\psi_1)\label{eq: proof c5-,1}\\
            &f_{12}+\psi_2=x_2(\partial_1\psi_1-\partial_2\psi_2)\label{eq: proof c5-,2}\\
            &f_{13}=(x_1-x_3)\partial_3\psi_1-x_2\partial_3\psi_2\label{eq: proof c5-,3}\\
            &f_{21}-\psi_2=x_2(\partial_1\psi_1-\partial_2\psi_2)\label{eq: proof c5-,4}\\
            &f_{22}-\psi_1=x_2(\partial_2\psi_1+\partial_1\psi_2)\label{eq: proof c5-,5}\\
            &f_{23}=(x_1-x_3)\partial_3\psi_2+x_2\partial_3\psi_1\label{eq: proof c5-,6}\\
            &f_{31}=(x_3-x_1)\partial_1\psi_3-x_2\partial_2\psi_3\label{eq: proof c5-,7}\\
            &f_{32}=(x_3-x_1)\partial_2\psi_3+x_2\partial_1\psi_3\label{eq: proof c5-,8}\\
            &f_{33}-\psi_3=0\label{eq: proof c5-,9}
        \end{empheq}
    \end{subequations}
    Equation \eqref{eq: proof c5-,9} gives 
    \begin{equation*}
        \psi_3=f_{33},
    \end{equation*}
    and taking the sum of \eqref{eq: proof c5-,1} and \eqref{eq: proof c5-,5} and the difference of \eqref{eq: proof c5-,2} and \eqref{eq: proof c5-,4}, we get respectively
    \begin{equation*}
        \psi_1=\frac{f_{11}+f_{22}}{2},\qquad \psi_2=\frac{f_{21}-f_{12}}{2},
    \end{equation*}
    i.e.\ we get the unique coordinate change as in the statement, and substituting these expressions in the equations \eqref{eq: conditions proof c5-} we get \eqref{eq: conditions for existence c5-}\footnote{Some equations turn out to give equivalent conditions, so we take the minimal set of independent conditions.}. 
    
\end{proof}

Now, it just remains to prove that \eqref{eq: conditions for existence c5-} just follows from the compatibility condition of $R_{\mathfrak{c}_5^-\oplus\mathfrak{d}}$ and $R_k$.

\begin{lemma}
    Let $[R_{\mathfrak{c}_5^-\oplus\mathfrak{d}},R_k]_{FN}=0$, then \eqref{eq: conditions for existence c5-} holds.
\end{lemma}

\begin{proof}
    Write for simplicity $[R_{\mathfrak{c}_5^-\oplus\mathfrak{d}},R_k]_{FN}=:B$, then vanishing of the Fr\"olicher-Nijenhuis brackets amounts to the following equations:
    \begin{align*}
        &B^1_{12} = x_2(\partial_2 f_{12}+\partial_1f_{11}-\partial_2f_{21}+\partial_1f_{22})-f_{21}-f_{12}=0\\
        &B^1_{13} = (x_1-x_3)\partial_3f_{11}+x_2(-\partial_3f_{21}+\partial_1f_{23}+\partial_2f_{13})-f_{13}=0\\
        &B^1_{23} = (x_1-x_3)\partial_3f_{12}+x_2(-\partial_3f_{22}+\partial_2f_{23}-\partial_1f_{13})+f_{23}=0\\
        &B^2_{12} = x_2(-\partial_1f_{12}+\partial_2f_{11}+\partial_2f_{22}+\partial_1 f_{21})+f_{11}-f_{22}=0\\
        &B^2_{13} = (x_1-x_3)\partial_3f_{21}+x_2(\partial_2f_{23}+\partial_3f_{11}-\partial_1f_{13})-f_{23}=0\\
        &B^2_{23} = (x_1-x_3)\partial_3f_{22}+x_2(-\partial_1f_{23}+\partial_3f_{12}-\partial_2f_{13})-f_{13}=0\\
        &B^3_{12} = (x_1-x_3)(\partial_1f_{32}-\partial_2f_{31})+x_2(\partial_1f_{31}+\partial_2f_{32})=0\\
        &B^3_{13} = (x_1-x_3)\partial_1f_{33}+x_2\partial_2f_{33}+f_{31}=0\\
        &B^3_{23} = (x_1-x_3)\partial_2f_{33}-x_2\partial_1f_{33}+f_{32}=0
    \end{align*}
    We see that 
    \begin{itemize}
        \item $B^1_{12}=0$ and $B^2_{12}=0$ are precisely \eqref{eq: c5-,2} and \eqref{eq: c5-,1};
        \item By taking the sum of $B^1_{13}=0$ and $B^2_{23}=0$ we get \eqref{eq: c5-,3};
        \item Taking the difference of $B^2_{13}=0$ and $B^1_{23}=0$ we get \eqref{eq: c5-,4};
        \item $B^3_{13}=0$ and $B^3_{23}=0$ are precisely \eqref{eq: c5-,5} and \eqref{eq: c5-,6}.
    \end{itemize}
\end{proof}

Summarizing, we have proven the following:

\begin{proposition}
    Let $\mathfrak{c}_5^-$ be the 2-dimensional LSA associated with linear operator
    \begin{equation*}
        R_{\mathfrak{c}_5^-}=\begin{pmatrix}
            x_1&-x_2\\x_2&x_1
        \end{pmatrix}
    \end{equation*}
    and let $\mathfrak{d}$ be a non-trivial one-dimensional LSA. Then $\mathfrak{c}_5^-\oplus\mathfrak{d}$ is non-degenerate in the formal category.
\end{proposition}

\subsection{Non-degeneracy of $\mathfrak{b}_4^+\oplus\mathfrak{d}$}

Let's consider the LSA $\mathfrak{b}_4^+\oplus\mathfrak{d}$, which is associated with the Nijenhuis operator 
\begin{equation}\label{eq: op b4+d}
    R_{\mathfrak{b}_4^+\oplus\mathfrak{d}}=\begin{pmatrix}
        2x_1&-x_2&0\\
        x_2&0&0\\
        0&0&x_3
    \end{pmatrix}
\end{equation}

\begin{lemma}
    The change of coordinates $y_i=x_i+\psi_i(x)$ with required condition \eqref{eq: condition lin} exists if and only if the components of $R_k$ satisfy the following relations:
     \begin{subequations}\label{eq: condition for existence b4+}
        \begin{empheq}[left=\empheqlbrace]{align}
            &x_2\mid f_{22}\label{eq: b4+,1}\\
            &2x_1\partial_1 f_{22}=-x_2(\partial_2f_{11}+\partial_2f_{22}+\partial_1f_{21}-\partial_1f_{12})\label{eq: b4+,2}\\
            &f_{21}+f_{12}=2x_1\partial_2f_{11}+x_2(\partial_1f_{11}+\partial_1f_{22}+\partial_2f_{12}-\partial_2f_{21})\label{eq: b4+,3}\\
            &2f_{13}=(2x_1-x_3)(\partial_3f_{11}+\partial_3f_{22})-x_2\left(\partial_3f_{21}-\partial_3f_{12}+2\dfrac{x_1}{x_2}\partial_3f_{22}\right)\label{eq: b4+,4}\\
            &2f_{23}=x_2(\partial_3f_{11}+\partial_3f_{22})-x_3\left(\partial_3f_{21}-\partial_3f_{12}+2\dfrac{x_1}{x_2}\partial_3f_{22}\right)\label{eq: b4+,5}\\
            &f_{31}=-(2x_1-x_3)\partial_1f_{33}-x_2\partial_2f_{33}\label{eq: b4+,6}\\
            &f_{32}=x_2\partial_1f_{33}+x_3\partial_2f_{33}\label{eq: b4+,7}
        \end{empheq}
    \end{subequations}
    In this case, the coordinate change takes the form
    \begin{equation*}
        \begin{cases}\vspace{0.1cm}
            y_1=x_1+\dfrac{f_{11}+f_{22}}{2}\\\vspace{0.1cm}
            y_2=x_2+\dfrac{f_{21}-f_{12}}{2}+x_1\dfrac{f_{22}}{x_2}\\
            y_3=x_3+f_{33}
        \end{cases}
    \end{equation*}
\end{lemma}

\begin{proof}
    For the operator \eqref{eq: op b4+d}, equation \eqref{eq: condition lin} is equivalent to
    \begin{subequations}\label{eq: condition proof b4+}
        \begin{empheq}[left=\empheqlbrace]{align}
            &f_{11}-2\psi_1=-x_2(\partial_1\psi_2+\partial_2\psi_1)\label{eq: proof b4+,1}\\
            &f_{12}+\psi_2=2x_1\partial_2\psi_1+x_2(\partial_1\psi_1-\partial_2\psi_2)\label{eq: proof b4+,2}\\
            &f_{13}=(2x_1-x_3)\partial_3\psi_1-x_2\partial_3\psi_2\label{eq: proof b4+,3}\\
            &f_{21}-\psi_2=x_2(\partial_1\psi_1-\partial_2\psi_2)-2x_1\partial_1\psi_2\label{eq: proof b4+,4}\\
            &f_{22}=x_2(\partial_2\psi_1+\partial_1\psi_2)\label{eq: proof b4+,5}\\
            &f_{23}=x_2\partial_3\psi_1-x_3\partial_3\psi_2\label{eq: proof b4+,6}\\
            &f_{31}=(x_3-2x_1)\partial_1\psi_3-x_2\partial_2\psi_3\label{eq: proof b4+,7}\\
            &f_{32}=x_3\partial_2\psi_3+x_2\partial_1\psi_3\label{eq: proof b4+,8}\\
            &f_{33}-\psi_3=0\label{eq: proof b4+,9}
        \end{empheq}
    \end{subequations}
    Equation \eqref{eq: proof b4+,9} gives 
    \begin{equation*}
        \psi_3=f_{33}.
    \end{equation*}
    If we consider the sum of \eqref{eq: proof b4+,1} and \eqref{eq: proof b4+,5}, we get
    \begin{equation*}
        \psi_1=\frac{f_{11}+f_{22}}{2}
    \end{equation*}
    Finally, taking the difference of \eqref{eq: proof b4+,2} and \eqref{eq: proof b4+,4}, we get
    \begin{equation*}
        2\psi_2+f_{12}-f_{21}=2x_1(\partial_2\psi_1+\partial_1\psi_2)
    \end{equation*}
    Using \eqref{eq: proof b4+,5} again, we get
    \begin{equation*}
        \partial_2\psi_1+\partial_1\psi_2=\frac{f_{22}}{x_2}
    \end{equation*}
    and hence substituting above we get
    \begin{equation*}
        2\psi_2=f_{21}-f_{12}+2\frac{x_1 f_{22}}{x_2}\quad\Longrightarrow\quad \psi_2=\frac{f_{21}-f_{12}}{2}+x_1\frac{f_{22}}{x_2}
    \end{equation*}
    Finally, substituting into the remaining equations \eqref{eq: condition proof b4+}, we get the relations in \eqref{eq: condition for existence b4+}, with $x_2\mid f_{22}$ to ensure that the coordinate change is well defined.
\end{proof}

Now, it just remains to prove that \eqref{eq: condition for existence b4+} just follows from the compatibility condition of $R_{\mathfrak{b}_4^+\oplus\mathfrak{d}}$ and $R_k$.

\begin{lemma}
    Let $[R_{\mathfrak{b}_4^+\oplus\mathfrak{d}},R_k]_{FN}=0$, then \eqref{eq: condition for existence b4+} holds.
\end{lemma}

\begin{proof}
     Write for simplicity $[R_{\mathfrak{b}_4^+\oplus\mathfrak{d}},R_k]_{FN}=:B$, then vanishing of the Fr\"olicher-Nijenhuis brackets amounts to the following equations:
     \begin{align*}
        &B^1_{12} = 2x_1\partial_2 f_{11}+x_2(\partial_2 f_{12}+\partial_1f_{11}-\partial_2f_{21}+\partial_1f_{22})-f_{21}-f_{12}=0\\
        &B^1_{13} = (2x_1-x_3)\partial_3f_{11}+x_2(-\partial_3f_{21}+\partial_1f_{23}+\partial_2f_{13})-2f_{13}=0\\
        &B^1_{23} = (2x_1-x_3)\partial_3f_{12}-2x_1\partial_2 f_{13}+x_2(-\partial_3f_{22}+\partial_2f_{23}-\partial_1f_{13})+f_{23}=0\\
        &B^2_{12} = x_2(-\partial_1f_{12}+\partial_2f_{11}+\partial_2f_{22}+\partial_1 f_{21})+2x_1\partial_1f_{22}=0\\
        &B^2_{13} = x_2(\partial_2f_{23}+\partial_3f_{11}-\partial_1f_{13})+2x_1\partial_1f_{23}-x_3\partial_3f_{21}-f_{23}=0\\
        &B^2_{23} = x_2(-\partial_2f_{13}+\partial_3f_{12}-\partial_1f_{23})-x_3\partial_3f_{22}=0\\
        &B^3_{12} = (2x_1-x_3)\partial_1f_{32}+x_2(\partial_1f_{31}+\partial_2f_{32})+x_3\partial_2 f_{31}+f_{32}=0\\
        &B^3_{13} = f_{31}+(2x_1-x_3)\partial_1f_{33}+x_2\partial_2f_{33}=0\\
        &B^3_{23} = f_{32}-x_2\partial_1f_{33}-x_3\partial_2f_{33}=0
    \end{align*}
    Now, we notice that
    \begin{itemize}
        \item using $B^2_{12}=0$ and $B^2_{23}=0$ we get
        \begin{align*}
            2x_1\partial_1f_{22}\equiv_{x_2}0,\qquad x_3\partial_3f_{22}\equiv_{x_2}0
        \end{align*}
        and hence
        \begin{equation*}
            x_2\mid\partial_1f_{22},\qquad x_2\mid\partial_3f_{22}
        \end{equation*}
        and it follows necessarily that $x_2\mid f_{22}$;
        \item $B^2_{12}=0$ is precisely \eqref{eq: b4+,2};
        \item $B^1_{12} =0$ is precisely \eqref{eq: b4+,3};
        \item $ B^3_{13}=0$ is precisely \eqref{eq: b4+,6};
        \item $B^3_{23}=0$ is precisely \eqref{eq: b4+,7};
        \item If we take the difference $B^2_{13}-B^1_{23}=0$ we get
    \begin{equation*}
        2f_{23}=x_2(\partial_3f_{11}+\partial_3f_{22})+2x_1(\partial_1f_{23}+\partial_2f_{13}-\partial_3f_{12})-x_3(\partial_3f_{21}-\partial_3f_{12}).
    \end{equation*}
    Thus, if we prove that
    \begin{equation*}
        -2\frac{x_1x_3}{x_2}\partial_3f_{22}=2x_1(\partial_1f_{23}+\partial_2f_{13}-\partial_3f_{12})
    \end{equation*}
    we get \eqref{eq: b4+,5}. But notice that such equality is equivalent to
    \begin{equation*}
        -x_3\partial_3f_{22}=x_2(\partial_1f_{23}+\partial_2f_{13}-\partial_3f_{12}),
    \end{equation*}
    and this is precisely $B^2_{23} = 0$;
    \item $B^1_{13}=0$ can be rewritten as
    \begin{equation*}
        2f_{13}=(2x_1-x_3)\partial_3f_{11}+x_2(\partial_1f_{23}-\partial_3f_{21}+\partial_2f_{13}),
    \end{equation*}
    thus it's enough to prove that
    \begin{equation*}
        (2x_1-x_3)\partial_3f_{22}+x_2\partial_3f_{12}-2x_1\partial_3f_{22}=x_2\partial_1f_{23}+x_2\partial_2f_{13}
    \end{equation*}
    to get the \eqref{eq: b4+,4}. But it's enough to use $B^2_{23}=0$, substituting the first term of the left hand side, and we are done.
    \end{itemize}
\end{proof}

Summarizing, we have proven the following:

\begin{proposition}
    Let $\mathfrak{b}_4^+$ be the 2-dimensional LSA associated with linear operator
    \begin{equation*}
        R_{\mathfrak{b}_4^+}=\begin{pmatrix}
            2x_1&-x_2\\x_2&0
        \end{pmatrix}
    \end{equation*}
    and let $\mathfrak{d}$ be a non-trivial one-dimensional LSA. Then $\mathfrak{b}_4^+\oplus\mathfrak{d}$ is non-degenerate in the formal category.
\end{proposition}

\subsection{Non-degeneracy of $\mathfrak{b}_4^-\oplus\mathfrak{d}$}

Let's now consider the LSA $\mathfrak{b}_4^-\oplus\mathfrak{d}$, which is associated with the Nijenhuis operator 
\begin{equation}\label{eq: op b4-d}
    R_{\mathfrak{b}_4^-\oplus\mathfrak{d}}=\begin{pmatrix}
        2x_1&x_2&0\\
        x_2&0&0\\
        0&0&x_3
    \end{pmatrix}
\end{equation}

\begin{lemma}
    The change of coordinates $y_i=x_i+\psi_i(x)$ with required condition \eqref{eq: condition lin} exists if and only if the components of $R_k$ satisfy the following relations:
    \begin{subequations}\label{eq: condition for existence b4-}
        \begin{empheq}[left=\empheqlbrace]{align}
            &x_2\mid f_{22}\label{eq: b4-,1}\\
            &2x_1\partial_1f_{22}=-x_2(\partial_2f_{11}+\partial_2f_{22}-\partial_{1}f_{12}-\partial_1f_{21})\label{eq: b4-,2}\\
            &f_{12}-f_{21}=2x_1\partial_2f_{11}-x_2(\partial_1f_{11}+\partial_1f_{22}-\partial_2f_{21}-\partial_2f_{12})\label{eq: b4-,3}\\
            &2f_{13}=2x_1\partial_3f_{11}+x_2\left(\partial_3f_{21}+\partial_3f_{12}\right)-x_3(\partial_3f_{11}+\partial_3f_{22})\label{eq: b4-,4}\\
            &2f_{23}=x_2(\partial_3f_{11}+\partial_3f_{22})-x_3\left(\partial_3f_{21}+\partial_3f_{12}-2\dfrac{x_1}{x_2}\partial_3f_{22}\right)\label{eq: b4-,5}\\
            &f_{31}=-(2x_1-x_3)\partial_1f_{33}-x_2\partial_2f_{33}\label{eq: b4-,6}\\
            &f_{32}=x_3\partial_2f_{33}-x_2\partial_1f_{33}\label{eq: b4-,7}
        \end{empheq}
    \end{subequations}
    In this case, the coordinate change takes the form
    \begin{equation*}
        \begin{cases}\vspace{0.1cm}
            y_1=x_1+\dfrac{f_{11}+f_{22}}{2}\\\vspace{0.1cm}
            y_2=x_2+\dfrac{f_{21}+f_{12}}{2}-x_1\dfrac{f_{22}}{x_2}\\
            y_3=x_3+f_{33}
        \end{cases}
    \end{equation*}
\end{lemma}

\begin{proof}
    For the operator \eqref{eq: op b4-d}, equation \eqref{eq: condition lin} is equivalent to
     \begin{subequations}\label{eq: condition proof b4-}
        \begin{empheq}[left=\empheqlbrace]{align}
            &f_{11}-2\psi_1=-x_2(\partial_2\psi_1-\partial_1\psi_2)\label{eq: proof b4-,1}\\
            &f_{12}-\psi_2=2x_1\partial_2\psi_1-x_2(\partial_1\psi_1-\partial_2\psi_2)\label{eq: proof b4-,2}\\
            &f_{13}=(2x_1-x_3)\partial_3\psi_1+x_2\partial_3\psi_2\label{eq: proof b4-,3}\\
            &f_{21}-\psi_2=x_2(\partial_1\psi_1-\partial_2\psi_2)-2x_1\partial_1\psi_2\label{eq: proof b4-,4}\\
            &f_{22}=x_2(\partial_2\psi_1-\partial_1\psi_2)\label{eq: proof b4-,5}\\
            &f_{23}=x_2\partial_3\psi_1-x_3\partial_3\psi_2\label{eq: proof b4-,6}\\
            &f_{31}=(x_3-2x_1)\partial_1\psi_3-x_2\partial_2\psi_3\label{eq: proof b4-,7}\\
            &f_{32}=x_3\partial_2\psi_3-x_2\partial_1\psi_3\label{eq: proof b4-,8}\\
            &f_{33}-\psi_3=0\label{eq: proof b4-,9}
        \end{empheq}
    \end{subequations}
    Equation \eqref{eq: proof b4-,9} gives 
    \begin{equation*}
        \psi_3=f_{33}.
    \end{equation*}
    If we consider the sum of \eqref{eq: proof b4-,1} and \eqref{eq: proof b4-,5}, we get
    \begin{equation*}
        \psi_1=\frac{f_{11}+f_{22}}{2}
    \end{equation*}
    Finally, taking the sum of \eqref{eq: proof b4-,2} and \eqref{eq: proof b4-,4}, we get
    \begin{equation*}
        -2\psi_2+f_{12}+f_{21}=2x_1(\partial_2\psi_1-\partial_1\psi_2)
    \end{equation*}
    Using \eqref{eq: proof b4-,5} again, we get
    \begin{equation*}
        \partial_2\psi_1-\partial_1\psi_2=\frac{f_{22}}{x_2}
    \end{equation*}
    and hence substituting above we get
    \begin{equation*}
        2\psi_2=f_{21}+f_{12}-2\frac{x_1 f_{22}}{x_2}\quad\Longrightarrow\quad \psi_2=\frac{f_{21}+f_{12}}{2}-x_1\frac{f_{22}}{x_2} 
    \end{equation*}
    Finally, substituting into the remaining equations \eqref{eq: condition proof b4-}, we get \eqref{eq: condition for existence b4-}, with $x_2\mid f_{22}$ to ensure that the coordinate change is well defined.
\end{proof}

Now, it just remains to prove that \eqref{eq: condition for existence b4-} just follows from the compatibility condition of $R_{\mathfrak{b}_4^-\oplus\mathfrak{d}}$ and $R_k$.

\begin{lemma}
    Let $[R_{\mathfrak{b}_4^-\oplus\mathfrak{d}},R_k]_{FN}=0$, then \eqref{eq: condition for existence b4-} holds.
\end{lemma}

\begin{proof}
    Write for simplicity $[R_{\mathfrak{b}_4^-\oplus\mathfrak{d}},R_k]_{FN}=:B$, then vanishing of the Fr\"olicher-Nijenhuis brackets amounts to the following equations:
    \begin{align*}
        &B^1_{12} = 2x_1\,\partial_2 f_{11} + x_2(-\partial_1 f_{11} - \partial_1 f_{22} + \partial_2 f_{12} + \partial_2 f_{21}) + f_{21} - f_{12} = 0 \\
        &B^1_{13} = (2x_1 - x_3)\,\partial_3 f_{11} + x_2(\partial_2 f_{13} + \partial_3 f_{21} - \partial_1 f_{23}) - 2f_{13} = 0 \\
        &B^1_{23} =(2x_1 - x_3)\,\partial_3 f_{12} - 2x_1\,\partial_2 f_{13} + x_2(\partial_1 f_{13} + \partial_3 f_{22} - \partial_2 f_{23}) - f_{23} = 0 \\
        &B^2_{12} =2x_1\,\partial_1 f_{22} + x_2(\partial_2 f_{11} - \partial_1 f_{12} - \partial_1 f_{21} + \partial_2 f_{22}) = 0 \\
        &B^2_{13} = 2x_1\,\partial_1 f_{23} + x_2(\partial_3 f_{11} - \partial_1 f_{13} + \partial_2 f_{23}) - x_3\,\partial_3 f_{21} - f_{23} = 0 \\
        &B^2_{23} = x_2(\partial_3 f_{12} - \partial_2 f_{13} + \partial_1 f_{23}) - x_3\,\partial_3 f_{22} = 0 \\
        &B^3_{12} =(2x_1 - x_3)\,\partial_1 f_{32} + x_2(\partial_2 f_{32} - \partial_1 f_{31}) + x_3\,\partial_2 f_{31} + f_{32} = 0 \\
        &B^3_{13} = (2x_1 - x_3)\,\partial_1 f_{33} + x_2\,\partial_2 f_{33} + f_{31} = 0 \\
        &B^3_{23} =x_2\,\partial_1 f_{33} - x_3\,\partial_2 f_{33} + f_{32} = 0
    \end{align*}
    Now, we notice that 
    \begin{itemize}
        \item Using $B^2_{12} =0$ and $B^2_{23}=0$ we get
    \begin{align*}
        2x_1\partial_1f_{22}\equiv_{x_2}0,\qquad x_3\partial_3f_{22}\equiv_{x_2}0
    \end{align*}
    and hence
    \begin{equation*}
        x_2\mid \partial_1f_{22},\qquad x_2\mid \partial_3f_{22}
    \end{equation*}
    and it follows necessarily that $x_2\mid f_{22}$;
    \item $B^2_{12}=0$ is precisely \eqref{eq: b4-,2};
    \item $B^1_{12} =0$ is precisely \eqref{eq: b4-,3};
    \item $B^3_{13}=0$ is precisely \eqref{eq: b4-,6};
    \item $B^3_{23} =0$ is \eqref{eq: b4-,7};
    \item If we take the sum of $B^1_{23}=0$ and $B^2_{13}=0$ we get
    \begin{equation*}
        2f_{23}=(2x_1-x_3)\partial_3f_{12}+2x_1(\partial_1f_{23}-\partial_2f_{13})+x_2(\partial_3f_{11}+\partial_3 f_{22})-x_3\partial_3f_{21}.
    \end{equation*}
    Thus, if we prove that
    \begin{equation*}
        (2x_1-x_3)\partial_3f_{12}+2x_1(\partial_1f_{23}-\partial_2f_{13})=-x_3\partial_3f_{12}+2\frac{x_1x_3}{x_2}\partial_3f_{22}
    \end{equation*}
    we get \eqref{eq: b4-,5}. But notice that such equality is equivalent to
    \begin{equation*}
        x_2(\partial_3f_{12}+\partial_1f_{23}-\partial_2f_{13})=x_3\partial_3f_{22}
    \end{equation*}
    and this is precisely $B^2_{23}=0$;
    \item Notice that \eqref{eq: b4-,4} is equivalent to $B^1_{13} + B^2_{23} =0$.
    \end{itemize}     
\end{proof}

Summarizing, we have proven the following:

\begin{proposition}
    Let $\mathfrak{b}_4^-$ be the 2-dimensional LSA associated with linear operator
    \begin{equation*}
        R_{\mathfrak{b}_4^-}=\begin{pmatrix}
            2x_1&x_2\\x_2&0
        \end{pmatrix}
    \end{equation*}
    and let $\mathfrak{d}$ be a non-trivial one-dimensional LSA. Then $\mathfrak{b}_4^-\oplus\mathfrak{d}$ is non-degenerate in the formal category.
\end{proposition}

\subsection{Non-degeneracy of $\mathfrak{s}_1$}

Let's consider the LSA $\mathfrak{s}_1$, which is associated with the Nijenhuis operator 
\begin{equation}\label{eq: op s1}
    R_{\mathfrak{s}_1}=\begin{pmatrix}
        2x_1-x_3&-x_2&x_3-x_1\\
        x_2&x_3&0\\
        0&0&x_3
    \end{pmatrix}
\end{equation}

\begin{lemma}
    The change of coordinates $y_i=x_i+\psi_i(x)$ with required condition \eqref{eq: condition lin} exists if and only if the components of $R_k$ satisfy the following relations:
    \begin{subequations}\label{eq: condition for existence s1}
        \begin{empheq}[left=\empheqlbrace]{align}
            &x_2\mid f_{32}\label{eq: s1,1}\\
            &x_2\mid f_{22}-f_{33}+\dfrac{f_{31}}{2}+2\dfrac{x_1-x_3}{x_2}f_{32}\label{eq: s1,2}\\
            &\begin{aligned}
                &2(x_1 - x_3)^2 \partial_1 f_{32} 
			+ x_2(x_1 - x_3) \left(2\partial_1 f_{22} 
			+  \partial_1 f_{31} +  \partial_2 f_{32}\right) 
			+ x_2^2 \left(\partial_2 f_{11} 
			- \partial_1 f_{12}
            + \partial_1 f_{21} 
            + \partial_2 f_{22}\right) \\
            & \qquad 
			+ (x_1 - x_3)f_{32} 
			+ x_2 f_{31}  = 0
            \end{aligned}\label{eq: s1,3}\\[2ex]
            &\begin{aligned}
            &(x_3 - x_1)\left( 2\partial_2 f_{11} 
		      - \partial_2 f_{31} 
		      + \partial_1 f_{32}\right) 
		      - x_2 \left( \partial_1 f_{11} 
		      + \partial_2 f_{12} 
              -  \partial_2 f_{21} 
		      +  \partial_1 f_{22} 
		      \right)\\
              &\qquad
		      + f_{12} 
		      + f_{21} 
		      - f_{32} = 0
              \end{aligned}\label{eq: s1,4}\\[2ex]
            &\begin{aligned}
            &-(x_1 - x_3)^2 \partial_1 f_{32} 
		      - x_2(x_1 - x_3) \left( \partial_1 f_{11} 
		      + 2 \partial_3 f_{11} 
		      +  \partial_1 f_{22} 
		      -  \partial_3 f_{31}\right) \\
              &\qquad
		      - x_2^2 \left( \partial_3 f_{12} 
		      - \partial_3 f_{21}\right) 
		      + 2(x_3 - x_1)f_{32} 
		      + x_2(f_{11} + 2f_{13} - f_{22} - f_{31})  = 0
              \end{aligned}\label{eq: s1,5}\\[2ex]
            &\begin{aligned}
            &- (x_1 - x_3)\bigg( 
		      \tfrac{1}{2}\left(\partial_1 f_{21} 
		      - \partial_1 f_{12}\right) 
		      + \tfrac{1}{x_2}\left(f_{22} - f_{33} + \tfrac{1}{2}f_{31} + \tfrac{2(x_1 - x_3)f_{32}}{x_2}\right) \\
              & \qquad
		      + \tfrac{x_1 - x_3}{x_2}\left(\partial_1 f_{22} - \partial_1 f_{33} + \tfrac{1}{2}\partial_1 f_{31} + \tfrac{2f_{32}}{x_2} + \tfrac{2(x_1 - x_3)\partial_1 f_{32}}{x_2}\right) \bigg) 
		      - \tfrac{1}{2}(x_1 - x_3)\partial_3 f_{32} \\
              &\qquad
		      - \tfrac{1}{2}x_2 \left(\partial_3 f_{11} + \partial_3 f_{22}\right) 
		      + \tfrac{1}{2}f_{32} +f_{23}= 0        
              \end{aligned}\label{eq: s1,6}\\[2ex]
            &-2(x_1 - x_3)^2 \partial_1 f_{32} 
		      - x_2(x_1 - x_3)\partial_2 f_{32} 
		      + 2x_2(x_1 - x_3)\partial_1 f_{33} 
		      + x_2^2 \partial_2 f_{33} 
		      + (x_3 - x_1)f_{32} 
		      + x_2 f_{31} = 0 
             \label{eq: s1,7} \\[2ex]
            &(x_1 - x_3)\partial_1 f_{32} - x_2 \partial_1 f_{33} + 2f_{32} = 0\label{eq: s1,8}
        \end{empheq}
    \end{subequations}
    In this case, the coordinate change takes the form
    \begin{equation*}
        \begin{cases}\vspace{0.1cm}
            y_1=x_1+\dfrac{1}{2}\left( f_{11}+ f_{22} + (x_1-x_3)\dfrac{f_{32}}{x_2}\right)\\\vspace{0.1cm}
            y_2=x_2+\dfrac{f_{21}-f_{12}}{2}+\dfrac{x_1-x_3}{x_2}\left(f_{22}-f_{33}+\dfrac{f_{31}}
            {2}+2(x_1-x_3)\dfrac{f_{32}}{x_2} \right)\\
            y_3=x_3+f_{33} -(x_1-x_3)\dfrac{f_{32}}{x_2}
        \end{cases}
    \end{equation*}
\end{lemma}

\begin{proof}
    For the operator \eqref{eq: op s1}, equation \eqref{eq: condition lin} is equivalent to
    \begin{subequations}\label{eq: conditions proof s1}
        \begin{empheq}[left=\empheqlbrace]{align}
            &f_{11} - 2 \, \psi_{1} + \psi_{3} = -(x_1-x_3){\partial_{1}\psi_{3}}-x_2 {\left({\partial_{1}\psi_{2}} + {\partial_{2}\psi_{1}}\right)} \label{eq: s11}\\
            &f_{12}+ \psi_{2} =(2x_1-2x_3)\partial_2\psi_1 +(x_3-x_1)\partial_2\psi_3+x_2(\partial_1\psi_1-\partial_2\psi_2) \label{eq: s12}\\
            &f_{13} + \psi_{1} - \psi_{3} = (x_1-x_3) \partial_1\psi_1+(2x_1-2x_3)\partial_3\psi_1+(x_3-x_1)\partial_3\psi_3-x_2\partial_3\psi_2 \label{eq: s13}\\
            &f_{21} - \psi_{2} = -2 \,(x_1-x_3) {\partial_{1}\psi_{2}} +x_2 {\left({\partial_{1}\psi_{1}} - {\partial_{2}\psi_{2}}\right)} \label{eq: s14}\\
            &f_{22} - \psi_{3} = x_2 {\left({\partial_{1}\psi_{2}} + {\partial_{2}\psi_{1}}\right)} \label{eq: s15}\\
            &f_{23} = (x_1-x_3){\partial_{1}\psi_{2}} + x_2{\partial_{3}\psi_{1}} \label{eq: s16}\\
            &f_{31} = -2 (x_1-x_3)\, {\partial_{1}\psi_{3}} - x_2 {\partial_{2}\psi_{3}} \label{eq: s17}\\
            &f_{32} = x_2{\partial_{1}\psi_{3}}  \label{eq: s18}\\
            &f_{33} - \psi_{3} = (x_1-x_3){\partial_{1}\psi_{3}} \label{eq: s19}
        \end{empheq}
    \end{subequations}
    Using \eqref{eq: s19} we get 
	\begin{equation*}
		\psi_3 = f_{33} - (x_1-x_3){\partial_{1}\psi_{3}} \overset{\eqref{eq: s18}}{=} f_{33} -\frac{x_1-x_3}{x_2} f_{32}.
	\end{equation*}
	Adding \eqref{eq: s11} and \eqref{eq: s15} we get
	\begin{equation*}
		\psi_1=\frac{1}{2}\left( f_{11}+ f_{22} + (x_1-x_3){\partial_{1}\psi_{3}}\right) \overset{\eqref{eq: s18}}{=} \frac{1}{2}\left( f_{11}+ f_{22} + \frac{x_1-x_3}{x_2} f_{32}\right).
	\end{equation*}
	Now, subtracting \eqref{eq: s14} from \eqref{eq: s12} we get
	\begin{equation*}
        2 \, \psi_2 = f_{21}-f_{12} + (2x_1-2x_3) \left( \, {\partial_{1}\psi_{2}} + {\partial_{2}\psi_{1}} \right) +(x_3-x_1)\partial_2\psi_3.
    \end{equation*}
	From \eqref{eq: s15} we have
    \begin{equation*}
        \partial_1\psi_2+\partial_2\psi_1=\frac{1}{x_2}(f_{22}-\psi_3)=\frac{1}{x_2}\left(f_{22}-f_{33}+\frac{x_1-x_3}{x_2}f_{32}\right)
    \end{equation*}
    and from \eqref{eq: s17} and using again \eqref{eq: s18} we get 
    \begin{equation*}
         -\frac{1}{x_2}\left(f_{31} + 2 \,\frac{x_1-x_3}{x_2} f_{32}\right) = {\partial_{2}\psi_{3}}.
    \end{equation*}
    Substituting these last expressions in the one above we get explicitly $\psi_2$: 
    \begin{equation*}
        \psi_2 = \frac{f_{21}-f_{12}}{2}+\frac{x_1-x_3}{x_2}\left(f_{22}-f_{33}+\frac{f_{31}}{2} +2\frac{x_1-x_3}{x_2}f_{32}\right).
    \end{equation*}
    Now, it's enough to insert the three expressions of $\psi_1,\psi_2,\psi_3$ in the equations \eqref{eq: conditions proof s1} to get the six compatibility conditions in \eqref{eq: condition for existence s1}, and the divisibility conditions are needed to ensure that the coordinate change is well defined.
\end{proof}

Now, it just remains to prove that \eqref{eq: condition for existence s1} just follows from the compatibility condition of $R_{\mathfrak{s}_1}$ and $R_k$.

\begin{lemma}
    Let $[R_{\mathfrak{s}_1},R_k]_{FN}=0$, then \eqref{eq: condition for existence s1} holds.
\end{lemma}

\begin{proof}
    Write for simplicity $[R_{\mathfrak{s}_1},R_k]_{FN}=:B$, then vanishing of the Fr\"olicher-Nijenhuis brackets amounts to the following equations:
    \begin{align*}
		B^1_{12} &= (2x_1 - 2x_3) \partial_2 f_{11} 
		+ (x_3 - x_1) \partial_2 f_{31} 
		+ (x_1 - x_3) \partial_1 f_{32} 
		+ x_2 \partial_1 f_{11} \\
		&\quad + x_2 \partial_2 f_{12} 
		- x_2 \partial_2 f_{21} 
		+ x_2 \partial_1 f_{22} 
		- f_{12} - f_{21} + f_{32}=0 \\[5pt]
		B^1_{13} &= (x_1 - x_3) \partial_1 f_{11} 
		+ (2x_1 - 2x_3) \partial_3 f_{11} 
		+ (x_3 - x_1) \partial_3 f_{31} 
		+ (x_1 - x_3) \partial_1 f_{33} \\
		&\quad + x_2 \partial_2 f_{13} 
		- x_2 \partial_3 f_{21} 
		+ x_2 \partial_1 f_{23} 
		- f_{11} - 2f_{13} + f_{31} + f_{33}=0 \\[5pt]
		B^1_{23} &= (x_1 - x_3) \partial_1 f_{12} 
		+ (2x_1 - 2x_3) \partial_3 f_{12} 
		+ (2x_3 - 2x_1) \partial_2 f_{13} 
		+ (x_3 - x_1) \partial_3 f_{32} \\
		&\quad + (x_1 - x_3) \partial_2 f_{33} 
		- x_2 \partial_1 f_{13} 
		- x_2 \partial_3 f_{22} 
		+ x_2 \partial_2 f_{23} 
		+ f_{23} + f_{32}=0\\[5pt]
		B^2_{12} &= (2x_1 - 2x_3) \partial_1 f_{22} 
		- x_2 \partial_1 f_{12} 
		+ x_2 \partial_1 f_{21} 
		+ x_2 \partial_2 f_{22} 
		+ x_2 \partial_2 f_{11} 
		+ f_{31}=0 \\[5pt]
		B^2_{13} &= (x_1 - x_3) \partial_1 f_{21} 
		+ (2x_1 - 2x_3) \partial_1 f_{23} 
		+ x_2 \partial_2 f_{23} 
		+ x_2 \partial_3 f_{11} 
		- x_2 \partial_1 f_{13} 
		- f_{23}=0 \\[5pt]
		B^2_{23} &= (x_1 - x_3) \partial_1 f_{22} 
		- x_2 \partial_2 f_{13} 
		- x_2 \partial_1 f_{23} 
		+ x_2 \partial_3 f_{12} 
		+ f_{22} - f_{33}=0\\[5pt]
		  B^3_{12} &= (2x_1 - 2x_3) \partial_1 f_{32} 
		+ x_2 \partial_1 f_{31} 
		+ x_2 \partial_2 f_{32} 
		+ f_{32}=0 \\[5pt]
		B^3_{13} &= (x_1 - x_3) \partial_1 f_{31} 
		+ (2x_1 - 2x_3) \partial_1 f_{33} 
		+ x_2 \partial_2 f_{33} 
		+ f_{31}=0 \\[5pt]
		B^3_{23} &= (x_1 - x_3) \partial_1 f_{32} 
		- x_2 \partial_1 f_{33} 
		+ 2f_{32}=0
	\end{align*}
    
    Then
    \begin{itemize}
        \item From $2 B^3_{23}-B^3_{12}$ we get
        \begin{equation*}
            f_{32}=\frac{x_2}{3}\left(2\partial_1 f_{33}+\partial_2 f_{32}+\partial_1 f_{31}\right),
        \end{equation*}
        i.e.\ $x_2\mid f_{32}$.
        \item Let us prove that 
        \begin{equation*}
        	x_2\mid f_{22}-f_{33}+\dfrac{f_{31}}{2}+2\dfrac{x_1-x_3}{x_2}f_{32}
        \end{equation*}
        This is the only non-trivial part of the proof. First, notice that from $B^2_{12}=0$ we can write
        \begin{equation}\label{eq: f_{31} s1}
        	f_{31}\equiv_{x_2}-2(x_1-x_3) \partial_1f_{22},
        \end{equation}
        and from $B^2_{23}=0$ we get
        \begin{equation}\label{eq: f_{22}-f_{33} s1}
        	f_{22}-f_{33}\equiv_{x_2}-(x_1-x_3)\partial_1f_{22}
        \end{equation}
        Thus
        \begin{equation*}
        	f_{22}-f_{33}+\frac{f_{31}}{2}\equiv_{x_2}-2(x_1-x_3)\partial_1f_{22}
        \end{equation*}
        We now use what we proved above, writing $f_{32}=x_2h$, so that $\partial_1f_{32}=x_2\partial_1h$ and $\partial_2f_{32}=h+x_2\partial_2h$, and from $B^3_{12}=0$ we get
        \begin{equation*}
        	2(x_1-x_3)\partial_1h+2h\equiv_{x_2}-\partial_1f_{31}
        \end{equation*}
        while from $B^3_{23}=0$ we get
        \begin{equation*}
        	2h+(x_1-x_3)\partial_1h=\partial_1f_{33}.
        \end{equation*}
        Now, taking the $\partial_1$ of the equations \eqref{eq: f_{31} s1} and \eqref{eq: f_{22}-f_{33} s1}, we get 
        \begin{align*}
        	\partial_1f_{31}&\equiv_{x_2}-2\partial_1f_{22}-2(x_1-x_3)\partial_1^2f_{22},\\
        	\partial_1f_{33}&\equiv_{x_2}2\partial_1f_{22}+(x_1-x_3)\partial_1^2f_{22}.
        \end{align*}
        Thus, substituting in the above we get the following equations:
        \begin{align*}
        	2(x_1-x_3)\partial_1h+2h&\equiv_{x_2}2\partial_1f_{22}+2(x_1-x_3)\partial_1^2f_{22},\\
        	(x_1-x_3)\partial_1h+2h&\equiv_{x_2}2\partial_1f_{22}+(x_1-x_3)\partial_1^2f_{22}.
        \end{align*}
        Thus, it follows immediately that
        \begin{equation*}
        	h\equiv_{x_2}\partial_1f_{22},
        \end{equation*}
        Thus, we have
        \begin{equation*}
        	2\frac{x_1-x_3}{x_2}f_{32}\equiv_{x_2}2(x_1-x_3) h=2(x_1-x_3) \partial_1f_{22}
        \end{equation*}
        and summarizing we have
        \begin{equation*}
        	f_{22}-f_{33}+\frac{1}{2}f_{31}+2\frac{x_1-x_3}{x_2}f_{32}\equiv_{x_2}0,
        \end{equation*}
        as claimed.
        \item Notice that \eqref{eq: s1,3} is equivalent to $(x_1-x_3)B^3_{12}+x_2 B^2_{12}=0$.
        \item Notice that \eqref{eq: s1,4} is equivalent to $-B^1_{12}=0$.
        \item Notice that \eqref{eq: s1,5} is equivalent to
        $x_2\left(-B^1_{13}-B^2_{23}-\frac{x_1-x_3}{x_2}B^3_{23}\right)=0.$
        \item Notice that \eqref{eq: s1,7} is equivalent to $B^3_{13}-\frac{x_1-x_3}{x_2}B^3_{12}=0.$
        \item Notice that \eqref{eq: s1,6} is equivalent to $\frac{1}{2}B^1_{23}-\frac{1}{2}B^2_{13}-\frac{x_1-x_3}{x_2}B^2_{23}-\frac{x_1-x_3}{2x_2}B^3_{13}-2\frac{(x_1-x_3)^2}{x_2^2}B^3_{23}=0.$
        \item Notice that \eqref{eq: s1,8} is equivalent to $B^3_{23}=0.$
    \end{itemize}
\end{proof}

Summarizing, we have proven the following:

\begin{proposition}
    Let $\mathfrak{s}_1$ be the 3-dimensional LSA associated with linear operator
    \begin{equation*}
        R_{\mathfrak{s}_1}=\begin{pmatrix}
        2x_1-x_3&-x_2&x_3-x_1\\
        x_2&x_3&0\\
        0&0&x_3
    \end{pmatrix}.
    \end{equation*}
    Then $\mathfrak{s}_1$ is non-degenerate in the formal category.
\end{proposition}

\subsection{Non-degeneracy of $\mathfrak{s}_2$}

Let's consider the LSA $\mathfrak{s}_2$, which is associated with the Nijenhuis operator 
\begin{equation}\label{eq: op s2}
    R_{\mathfrak{s}_2}=\begin{pmatrix}
        2x_1-x_3&x_2&x_3-x_1\\
        x_2&x_3&0\\
        0&0&x_3
    \end{pmatrix}
\end{equation}

\begin{lemma}
    The change of coordinates $y_i=x_i+\psi_i(x)$ with required condition \eqref{eq: condition lin} exists if and only if the components of $R_k$ satisfy the following relations:
    \begin{subequations}\label{eq: condition for existence s2}
        \begin{empheq}[left=\empheqlbrace]{align}
            &x_2\mid f_{32}\label{eq: s2,1}\\
            &x_2\mid f_{22}-f_{33}+\dfrac{1}{2}f_{31}+2\dfrac{x_3-x_1}{x_2}f_{32}\label{eq: s2,2}\\
            &\begin{aligned}
                &-2(x_1 - x_3)^2 \partial_1 f_{32} 
            + 2x_2(x_1 - x_3) \partial_1 f_{22} 
            + x_2(x_1 - x_3) \partial_1 f_{31} 
            - x_2(x_1 - x_3) \partial_2 f_{32} 
            + x_2^2 \partial_2 f_{11} 
            - x_2^2 \partial_1 f_{12} \\
            &\qquad
            - x_2^2 \partial_1 f_{21} 
            + x_2^2 \partial_2 f_{22} 
            + (x_3 - x_1)f_{32} 
            + x_2 f_{31} = 0 \end{aligned}\label{eq: s2,3}\\[2ex]
            &\begin{aligned}
            &2(x_3 - x_1) \partial_2 f_{11} 
            + (x_1 - x_3) \partial_2 f_{31} 
            + (x_3 - x_1) \partial_1 f_{32} 
            + x_2 \partial_1 f_{11} 
            - x_2 \partial_2 f_{12} \\
            &\qquad
            - x_2 \partial_2 f_{21} 
            + x_2 \partial_1 f_{22}
            + f_{12} 
            - f_{21} 
            - f_{32} = 0 \end{aligned}\label{eq: s2,4}\\[2ex]
            &\begin{aligned}
            &(x_1 - x_3)^2 \partial_1 f_{32} 
            - x_2(x_1 - x_3) \partial_1 f_{11} 
            - 2x_2(x_1 - x_3) \partial_3 f_{11} 
            - x_2(x_1 - x_3) \partial_1 f_{22} 
            + x_2(x_1 - x_3) \partial_3 f_{31} \\
            &\qquad
            - x_2^2 \partial_3 f_{12} 
            - x_2^2 \partial_3 f_{21} 
            + 2(x_1 - x_3)f_{32} 
            + x_2(f_{11} + 2f_{13} - f_{22} - f_{31})= 0\end{aligned} \label{eq: s2,5}\\[2ex]
            &\begin{aligned}
            &2f_{23} - 2(x_1 - x_3) \Bigg[ 
            \dfrac{1}{2}\partial_1 f_{21} 
            + \dfrac{1}{2}\partial_1 f_{12} 
            + \dfrac{1}{x_2}\left(f_{33} - f_{22} - \dfrac{1}{2}f_{31}\right) 
            + \dfrac{x_1 - x_3}{x_2}\left(\partial_1 f_{33} - \partial_1 f_{22} - \frac{1}{2}\partial_1 f_{31}\right)\\
            &\qquad
            + \dfrac{4(x_1 - x_3)f_{32}}{x_2^2} 
            + \dfrac{2(x_1 - x_3)^2 \partial_1 f_{32}}{x_2^2} \Bigg] 
            - x_2 \left[ \partial_3 f_{11} + \partial_3 f_{22} + \dfrac{f_{32}}{x_2} +
            \frac{(x_3 - x_1)\partial_3 f_{32}}{x_2} \right] = 0\end{aligned}\label{eq: s2,6}\\[2ex]
            &2(x_1 - x_3)^2 \partial_1 f_{32} 
            + x_2(x_1 - x_3) \partial_2 f_{32} 
            + 2x_2(x_1 - x_3) \partial_1 f_{33} 
            + x_2^2 \partial_2 f_{33} 
            + (x_1 - x_3)f_{32} 
            + x_2 f_{31} = 0\label{eq: s2,7}\\[2ex]
            &(x_1 - x_3) \partial_1 f_{32} 
            + x_2 \partial_1 f_{33} 
            + 2f_{32} = 0 \label{eq: s2,8}
        \end{empheq}
    \end{subequations}
    In this case, the coordinate change takes the form
    \begin{equation*}
        \begin{cases}\vspace{0.1cm}
            y_1=x_1+\dfrac{1}{2}\left( f_{11}+ f_{22} + \dfrac{x_3-x_1}{x_2} f_{32}\right)\\\vspace{0.1cm}
            y_2=x_2+\dfrac{f_{12}+f_{21}}{2}+\dfrac{x_3-x_1}{x_2}\left(f_{22}-f_{33}+\dfrac{1}{2}f_{31}+2\dfrac{x_3-x_1}{x_2}f_{32}\right)\\
            y_3=x_3+f_{33}+\dfrac{x_1-x_3}{x_2} f_{32}
        \end{cases}
    \end{equation*}
\end{lemma}

\begin{proof}
    For the operator \eqref{eq: op s2}, equation \eqref{eq: condition lin} is equivalent to
    \begin{subequations}\label{eq: conditions proof s2}
        \begin{empheq}[left=\empheqlbrace]{align}
        &f_{11} - 2 \, \psi_{1} + \psi_{3} = -{\partial_{1}\psi_{3}} x_{1} + {\left({\partial_{1}\psi_{2}} - {\partial_{2}\psi_{1}}\right)} x_{2} + {\partial_{1}\psi_{3}} x_{3} \label{eq: s21}\\
        &f_{12} - \psi_{2} = 2 \, {\partial_{2}\psi_{1}} x_{1} - {\partial_{2}\psi_{3}} x_{1} - {\left({\partial_{1}\psi_{1}} - {\partial_{2}\psi_{2}}\right)} x_{2} - {\left(2 \, {\partial_{2}\psi_{1}} - {\partial_{2}\psi_{3}}\right)} x_{3}\label{eq: s22}\\
        & f_{13} + \psi_{1} - \psi_{3} = {\partial_{1}\psi_{1}} x_{1} + 2 \, {\partial_{3}\psi_{1}} x_{1} - {\partial_{3}\psi_{3}} x_{1} + {\partial_{3}\psi_{2}} x_{2} - {\left({\partial_{1}\psi_{1}} + 2 \, {\partial_{3}\psi_{1}} - {\partial_{3}\psi_{3}}\right)} x_{3} \label{eq: s23}\\
        &f_{21} - \psi_{2} = -2 \, {\partial_{1}\psi_{2}} x_{1} + {\left({\partial_{1}\psi_{1}} - {\partial_{2}\psi_{2}}\right)} x_{2} + 2 \, {\partial_{1}\psi_{2}} x_{3} \label{eq: s24}\\
        &f_{22} - \psi_{3} = -{\left({\partial_{1}\psi_{2}} - {\partial_{2}\psi_{1}}\right)} x_{2} \label{eq: s25}\\
        &f_{23} = {\partial_{1}\psi_{2}} x_{1} + {\partial_{3}\psi_{1}} x_{2} - {\partial_{1}\psi_{2}} x_{3}\label{eq: s26}\\
        &f_{31} = -2 \, {\partial_{1}\psi_{3}} x_{1} - {\partial_{2}\psi_{3}} x_{2} + 2 \, {\partial_{1}\psi_{3}} x_{3} \label{eq: s27}\\
        &f_{32} = -{\partial_{1}\psi_{3}} x_{2}\label{eq: s28}\\	
        &f_{33} - \psi_{3} = {\partial_{1}\psi_{3}} x_{1} - {\partial_{1}\psi_{3}} x_{3} \label{eq: s29}
        \end{empheq}
    \end{subequations}
    Using \eqref{eq: s29} we get 
	\begin{equation*}
		\psi_3 = f_{33} - (x_1-x_3){\partial_{1}\psi_{3}} \overset{\eqref{eq: s28}}{=} f_{33}+\frac{x_1-x_3}{x_2} f_{32}.
	\end{equation*}
    Adding \eqref{eq: s21} and \eqref{eq: s25} we get
	\begin{equation*}
		2\psi_1=f_{11}+ f_{22} + (x_1-x_3){\partial_{1}\psi_{3}} \overset{\eqref{eq: s28}}{=}  f_{11}+ f_{22} + \frac{x_3-x_1}{x_2} f_{32}.
	\end{equation*}
    Now, using equations \eqref{eq: s24}, \eqref{eq: s22}, \eqref{eq: s25} and \eqref{eq: s27} we get
    \begin{equation*}
        \psi_2=\frac{f_{12}+f_{21}}{2}+\frac{x_3-x_1}{x_2}\left(f_{22}-f_{33}+\dfrac{1}{2}f_{31}+2\dfrac{x_3-x_1}{x_2}f_{32}\right).
    \end{equation*}
    Now, it's enough to insert the three expressions of $\psi_1,\psi_2,\psi_3$ in the equations \eqref{eq: conditions proof s2} to get the six compatibility conditions in \eqref{eq: condition for existence s2}, with the divisibility conditions needed to ensure that the coordinate change is well defined.
\end{proof}

Now, it just remains to prove that \eqref{eq: condition for existence s2} just follows from the compatibility condition of $R_{\mathfrak{s}_2}$ and $R_k$.

\begin{lemma}
    Let $[R_{\mathfrak{s}_2},R_k]_{FN}=0$, then \eqref{eq: condition for existence s2} holds.
\end{lemma}

\begin{proof}
    Write for simplicity $[R_{\mathfrak{s}_2},R_k]_{FN}=:B$, then vanishing of the Fr\"olicher-Nijenhuis brackets amounts to the following equations:
    \begin{align*}
		B^1_{12} &= 2(x_1 - x_3) \partial_2 f_{11} 
    + (x_3 - x_1) \partial_2 f_{31} 
    + (x_1 - x_3) \partial_1 f_{32} \\
    &- x_2 \partial_1 f_{11} 
    + x_2 \partial_2 f_{12} 
    + x_2 \partial_2 f_{21} 
    - x_2 \partial_1 f_{22} 
    - f_{12} 
    + f_{21} 
    + f_{32} =0 \\[5pt]
		B^1_{13} &= (x_1 - x_3) \partial_1 f_{11} 
    + 2(x_1 - x_3) \partial_3 f_{11} 
    + (x_3 - x_1) \partial_3 f_{31} 
    + (x_1 - x_3) \partial_1 f_{33} \\
    &+ x_2 \partial_2 f_{13} 
    + x_2 \partial_3 f_{21} 
    - x_2 \partial_1 f_{23} 
    - f_{11} 
    - 2f_{13} 
    + f_{31} 
    + f_{33}=0 \\[5pt]
		B^1_{23} &= (x_1 - x_3) \partial_1 f_{12} 
    + 2(x_1 - x_3) \partial_3 f_{12} 
    - 2(x_1 - x_3) \partial_2 f_{13} 
    + (x_3 - x_1) \partial_3 f_{32} \\
    &+ (x_1 - x_3) \partial_2 f_{33} 
    + x_2 \partial_1 f_{13} 
    + x_2 \partial_3 f_{22} 
    - x_2 \partial_2 f_{23} 
    - f_{23} 
    + f_{32}=0\\[5pt]
		B^2_{12} &= 2(x_1 - x_3) \partial_1 f_{22} 
    - x_2 \partial_1 f_{12} 
    - x_2 \partial_1 f_{21} 
    + x_2 \partial_2 f_{22} 
    + x_2 \partial_2 f_{11} 
    + f_{31}=0 
        \end{align*}
        \begin{align*}
		B^2_{13} &= (x_1 - x_3) \partial_1 f_{21} 
    + 2(x_1 - x_3) \partial_1 f_{23} 
    + x_2 \partial_2 f_{23} 
    + x_2 \partial_3 f_{11} 
    - x_2 \partial_1 f_{13} 
    - f_{23}=0 \\[5pt]
		B^2_{23} &= (x_1 - x_3) \partial_1 f_{22} 
    - x_2 \partial_2 f_{13} 
    + x_2 \partial_1 f_{23}
    + x_2 \partial_3 f_{12} 
    + f_{22} 
    - f_{33}=0\\[5pt]
		  B^3_{12} &= 2(x_1 - x_3) \partial_1 f_{32} 
    - x_2 \partial_1 f_{31} 
    + x_2 \partial_2 f_{32} 
    + f_{32}=0 \\[5pt]
		B^3_{13} &= (x_1 - x_3) \partial_1 f_{31} 
    + 2(x_1 - x_3) \partial_1 f_{33} 
    + x_2 \partial_2 f_{33} 
    + f_{31}=0 \\[5pt]
		B^3_{23} &= (x_1 - x_3) \partial_1 f_{32} 
    + x_2 \partial_1 f_{33} 
    + 2f_{32}=0
	\end{align*}
    Then
    \begin{itemize}
        \item Computing $2B^3_{23}-B^3_{12}$ we get
        \begin{equation*}
            3f_{32}=-x_2\left(2\partial_1f_{33}+\partial_1f_{31}-\partial_2f_{32}\right)
        \end{equation*}
        and hence $x_2\mid f_{32}$.
        \item Let us prove that
        \begin{equation*}
        	x_2\mid f_{22}-f_{33}+\dfrac{1}{2}f_{31}+2\dfrac{x_3-x_1}{x_2}f_{32}
        \end{equation*}
        From the system, it's easy to show, using $B^2_{23}=0$ and $B^2_{12}=0$ that
        \begin{equation}
        	f_{22}-f_{33}+\frac{1}{2}f_{31}\equiv_{x_2}2(x_3-x_1)\partial_1f_{22}
        \end{equation}
        Now, if we write $f_{32}=x_2h$ then from $B^3_{23}=0$ we get
        \begin{equation*}
        	2h=(x_3-x_1)\partial_1h-\partial_1f_{33},
        \end{equation*}
        while from $B^3_{12}=0$ we get
        \begin{equation*}
        	2h\equiv_{x_2}2(x_3-x_1)\partial_1h+\partial_1f_{31}
        \end{equation*}
        and hence, it follows that
        \begin{equation*}
        	2h\equiv_{x_2}-(2\partial_1f_{33}+\partial_1f_{31}).
        \end{equation*}
        From $\partial_1B^2_{12}=0$ and $\partial_1 B^2_{23}=0$ we get
        \begin{equation*}
        	\partial_1f_{31}\equiv_{x_2}-2\partial_1f_{22}+2(x_3-x_1)\partial_1^2f_{22}
        \end{equation*}
        and
        \begin{equation*}
        	\partial_1f_{33}\equiv_{x_2}2\partial_1f_{22}-(x_3-x_1)\partial_1^2f_{22}
        \end{equation*}
        Hence, it follows that
        \begin{equation*}
        	2\partial_1f_{33}+\partial_1f_{31}\equiv_{x_2}2\partial_1f_{22}
        \end{equation*}
        and hence
        \begin{equation*}
        	2h\equiv_{x_2}-2\partial_1f_{22}.
        \end{equation*}
        Finally, we can conclude that
        \begin{align*}
        	f_{22}-f_{33}+\frac{1}{2}f_{31}+2\frac{x_3-x_1}{x_2}f_{32}\equiv_{x_2}2(x_3-x_1)\partial_1 f_{22}-2(x_3-x_1)\partial_1f_{22}=0.
        \end{align*}
        \item Notice that \eqref{eq: s2,3} is equivalent to $(x_3-x_1)B^3_{12}+x_2 B^2_{12}=0.$ 
        \item Notice that \eqref{eq: s2,4} is equivalent to $-B^1_{12}=0.$
        \item Notice that \eqref{eq: s2,5} is equivalent to $x_2\left(-B^1_{13}-B^2_{23}+\frac{x_1-x_3}{x_2}B^3_{23}\right)=0.$
        \item Notice that \eqref{eq: s2,6} is equivalent to $-2\bigg[\frac{1}{2}B^1_{23}+\frac{1}{2}B^2_{13}+\frac{x_3-x_1}{x_2}B^2_{23}+\frac{x_3-x_1}{2x_2}B^3_{13}+2\frac{(x_1-x_3)^2}{x_2^2}B^3_{23}\bigg]=0.$
        \item Notice that \eqref{eq: s2,7} is equivalent to $x_2(B^3_{13}+\frac{x_1-x_3}{x_2}B^3_{12})=0.$
        \item Notice that \eqref{eq: s2,8} is equivalent to $B^3_{23}=0.$
    \end{itemize}
\end{proof}

Summarizing, we have proven the following:

\begin{proposition}
    Let $\mathfrak{s}_2$ be the 3-dimensional LSA associated with linear operator
    \begin{equation*}
        R_{\mathfrak{s}_2}=\begin{pmatrix}
        2x_1-x_3&x_2&x_3-x_1\\
        x_2&x_3&0\\
        0&0&x_3
    \end{pmatrix}.
    \end{equation*}
    Then $\mathfrak{s}_2$ is non-degenerate in the formal category.
\end{proposition}

\subsection{Non-degeneracy of $\mathfrak{s}_3$}

Let's consider the LSA $\mathfrak{s}_3$, which is associated with the Nijenhuis operator 
\begin{equation}\label{eq: op s3}
    R_{\mathfrak{s}_3}=\begin{pmatrix}
        x_1&-x_3&-x_2\\
        x_2&0&-x_2\\
        0&0&x_3
    \end{pmatrix}
\end{equation}

\begin{lemma}
    The change of coordinates $y_i=x_i+\psi_i(x)$ with required condition \eqref{eq: condition lin} exists if and only if the components of $R_k$ satisfy the following relations:
    \begin{subequations}\label{eq: condition for existence s3}
        \begin{empheq}[left=\empheqlbrace]{align}
            & x_3\mid f_{32} \label{eq: cond s3,1}\\
            & x_3\mid x_2(x_1+x_2)\frac{f_{32}}{x_{3}} + x_1 f_{22} - x_2 f_{12} + x_2 f_{31} - x_2 f_{33} \label{eq: cond s3,2}\\
            & \partial_{2}f_{11} x_{1} + \partial_{1}f_{32} x_{2} + \partial_{2}f_{12} x_{2} - \partial_{2}f_{31} x_{2} + \partial_{1}f_{11} x_{3} + \partial_{1}f_{22} x_{3} - \partial_{2}f_{21} x_{3} - f_{31}=0 \label{eq: cond s3,3}\\[2ex]
            & \begin{aligned}
                &\partial_{1}f_{32} x_{2}^{2} + \partial_{2}f_{32} x_{2}^{2} + \partial_{3}f_{11} x_{1} x_{3} + \partial_{1}f_{11} x_{2} x_{3} + \partial_{1}f_{22} x_{2} x_{3} + \partial_{2}f_{11} x_{2} x_{3} + \partial_{2}f_{22} x_{2} x_{3} \\
                &\qquad
                + \partial_{3}f_{12} x_{2} x_{3} - \partial_{3}f_{31} x_{2} x_{3} - \partial_{3}f_{32} x_{2} x_{3} - \partial_{3}f_{11} x_{3}^{2} - \partial_{3}f_{21} x_{3}^{2} - \partial_{3}f_{22} x_{3}^{2} + 2 \, f_{32} x_{2} - f_{13} x_{3} - f_{21} x_{3}=0 
            \end{aligned} \label{eq: cond s3,4}\\[2ex]
            & \begin{aligned}
                &\partial_{1}f_{32} x_{1} x_{2} + \partial_{1}f_{32} x_{2}^{2} + 2 \, \partial_{2}f_{32} x_{2}^{2} + \partial_{1}f_{22} x_{1} x_{3} - \partial_{1}f_{12} x_{2} x_{3} + \partial_{1}f_{31} x_{2} x_{3} - \partial_{1}f_{33} x_{2} x_{3} \\
                &\qquad
                + \partial_{2}f_{11} x_{2} x_{3} + \partial_{2}f_{22} x_{2} x_{3} - \partial_{2}f_{33} x_{2} x_{3} + \partial_{1}f_{21} x_{3}^{2} + 3 \, f_{32} x_{2}=0 
            \end{aligned} \label{eq: cond s3,5}\\[2ex]
            & \begin{aligned}
                &\partial_{1}f_{32} x_{1} x_{2}^{2} + \partial_{2}f_{32} x_{1} x_{2}^{2} + \partial_{1}f_{32} x_{2}^{3} + \partial_{2}f_{32} x_{2}^{3} + \partial_{1}f_{22} x_{1} x_{2} x_{3} + \partial_{2}f_{22} x_{1} x_{2} x_{3} - \partial_{3}f_{32} x_{1} x_{2} x_{3} \\
                &\qquad
                - \partial_{1}f_{12} x_{2}^{2} x_{3} + \partial_{1}f_{31} x_{2}^{2} x_{3} - \partial_{1}f_{33} x_{2}^{2} x_{3} - \partial_{2}f_{12} x_{2}^{2} x_{3} + \partial_{2}f_{31} x_{2}^{2} x_{3} - \partial_{2}f_{33} x_{2}^{2} x_{3} \\
                &\qquad
                + \partial_{3}f_{32} x_{2}^{2} x_{3} - \partial_{3}f_{22} x_{1} x_{3}^{2} + \partial_{1}f_{21} x_{2} x_{3}^{2} + \partial_{2}f_{21} x_{2} x_{3}^{2} + \partial_{3}f_{11} x_{2} x_{3}^{2} + \partial_{3}f_{12} x_{2} x_{3}^{2} \\
                &\qquad
                + \partial_{3}f_{22} x_{2} x_{3}^{2} - \partial_{3}f_{31} x_{2} x_{3}^{2} - \partial_{3}f_{21} x_{3}^{3} + 2 \, f_{32} x_{1} x_{2} + 2 \, f_{32} x_{2}^{2} - f_{12} x_{2} x_{3} + f_{22} x_{2} x_{3} \\
                &\qquad
                + f_{31} x_{2} x_{3} - f_{33} x_{2} x_{3} - f_{21} x_{3}^{2} - f_{23} x_{3}^{2}=0 
            \end{aligned} \label{eq: cond s3,6}\\[2ex]
            & \partial_{1}f_{32} x_{1} x_{2} + \partial_{2}f_{32} x_{2}^{2} - \partial_{1}f_{33} x_{1} x_{3} - \partial_{1}f_{32} x_{2} x_{3} - \partial_{2}f_{33} x_{2} x_{3} + \partial_{1}f_{33} x_{3}^{2} + f_{32} x_{2} - f_{31} x_{3}=0 \label{eq: cond s3,7}\\[2ex]
            & -\partial_{1}f_{32} x_{2} - \partial_{2}f_{32} x_{2} + \partial_{1}f_{33} x_{3} + \partial_{2}f_{33} x_{3} - 2 \, f_{32}=0\label{eq: cond s3,8}
        \end{empheq}
    \end{subequations}
    In this case, the coordinate change takes the form
    \begin{equation*}
        \begin{cases}\vspace{0.1cm}
            y_1=x_1+f_{11} + f_{22} + x_2 \dfrac{f_{32}}{x_3}\\\vspace{0.1cm}
            y_2=x_2+f_{21} + \dfrac{1}{x_3}\left(x_2(x_1+x_2)\dfrac{f_{32}}{x_{3}} +  x_1 f_{22} - x_2 f_{12} + x_2 f_{31} - x_2 f_{33}\right)\\
            y_3=x_3+f_{33}-x_2 \dfrac{f_{32}}{x_3}
        \end{cases}
    \end{equation*}
\end{lemma}

\begin{proof}
    For the operator \eqref{eq: op s3}, equation \eqref{eq: condition lin} is equivalent to
    \begin{subequations}\label{eq: conditions proof s3}
        \begin{empheq}[left=\empheqlbrace]{align}
		f_{11} - \psi_1 &= -x_3 \partial_1 \psi_2 - x_2 \left( \partial_1 \psi_3 + \partial_2 \psi_1 \right) \label{eq: s31} \\
		f_{12} + \psi_3 &= x_1 \partial_2 \psi_1 - x_2 \partial_2 \psi_3 + x_3 \left( \partial_1 \psi_1 - \partial_2 \psi_2 \right) \label{eq: s32} \\
		f_{13} + \psi_2 &= (x_1 - x_3)\partial_3 \psi_1 + x_2 \left( \partial_1 \psi_1 + \partial_2 \psi_1 - \partial_3 \psi_3 \right) - x_3 \partial_3 \psi_2 \label{eq: s33} \\
		f_{21} - \psi_2 &= -x_1 \partial_1 \psi_2 + x_2 \left( \partial_1 \psi_1 - \partial_2 \psi_2 - \partial_1 \psi_3 \right) \label{eq: s34} \\
		f_{22} &= x_2 \left( \partial_2 \psi_1 - \partial_2 \psi_3 \right) + x_3 \partial_1 \psi_2 \label{eq: s35} \\
		f_{23} + \psi_2 &= x_2 \left( \partial_3 \psi_1 - \partial_3 \psi_3 + \partial_1 \psi_2 + \partial_2 \psi_2 \right) - x_3 \partial_3 \psi_2 \label{eq: s36} \\
		f_{31} &= (x_3 - x_1)\partial_1 \psi_3 - x_2 \partial_2 \psi_3 \label{eq: s37} \\
		f_{32} &= x_3 \left( \partial_1 \psi_3 + \partial_2 \psi_3 \right) \label{eq: s38} \\
		f_{33} - \psi_3 &= x_2 \left( \partial_1 \psi_3 + \partial_2 \psi_3 \right) \label{eq: s39}
        \end{empheq}
    \end{subequations}
	From \eqref{eq: s38}
	\begin{equation*}
		x_3|f_{32},\qquad\frac{f_{32}}{x_3} = \partial_1 \psi_3 + \partial_2 \psi_3
	\end{equation*}
	and from \eqref{eq: s39}
	\begin{align*}
		\psi_3 = f_{33} - x_2 \left( \partial_1 \psi_3 + \partial_2 \psi_3 \right)\quad \Longrightarrow \quad \psi_3 = f_{33} - \frac{x_2}{x_3} f_{32}
	\end{align*}
	and we get the explicit expression of $\psi_3$, while using \eqref{eq: s31} we get
	\begin{align*}
		\psi_1 &= f_{11} + x_3 \partial_1 \psi_2 + x_2 \left( \partial_1 \psi_3 + \partial_2 \psi_1 \right) \\
		&\overset{\eqref{eq: s35}}{=} f_{11} + f_{22} - x_2 \left( \partial_2 \psi_1 - \partial_2 \psi_3 \right) + x_2 \left( \partial_1 \psi_3 + \partial_2 \psi_1 \right) \\
		&=f_{11}+f_{22}+x_2(\partial_1\psi_3+\partial_2\psi_3)
	\end{align*}
	and using \eqref{eq: s38} we get
	\begin{align*}
		{\psi_1} &{= f_{11} + f_{22} + x_2 \frac{f_{32}}{x_3}}
	\end{align*}
	and we get the explicit expression of $\psi_1$. Now, from \eqref{eq: s34} we get
    \begin{equation*}
        \psi_2=f_{21}+x_1\partial_1\psi_2-x_2(\partial_1\psi_1-\partial_2\psi_2-\partial_1\psi_3).
    \end{equation*}
    We now claim that
    \begin{equation}\label{eq: cond for psi2}
        x_1\partial_1\psi_2-x_2(\partial_1\psi_1-\partial_2\psi_2-\partial_1\psi_3)=\dfrac{1}{x_3}\left(x_2(x_1+x_2)\dfrac{f_{32}}{x_{3}} +  x_1 f_{22} - x_2 f_{12} + x_2 f_{31} - x_2 f_{33}\right)
    \end{equation}
    To prove that, notice that
    \begin{enumerate}
        \item From \eqref{eq: s35} we get
        \begin{equation*}
            x_1f_{22}=x_1x_2(\partial_2\psi_1-\partial_2\psi_3)+x_1x_3\partial_1\psi_2;
        \end{equation*}
        \item From \eqref{eq: s32} we get
        \begin{equation*}
            x_2f_{12}=x_1x_2\partial_2\psi_1-x_2\psi_3-x_2^2\partial_2\psi_3+x_2x_3(\partial_1\psi_1-\partial_2\psi_2);
        \end{equation*}
        \item From \eqref{eq: s37} we get
        \begin{equation*}
            x_2f_{31}=x_2(x_3-x_1)\partial_1\psi_3-x_2^2\partial_2\psi_3;
        \end{equation*}
        \item From \eqref{eq: s39} we get
        \begin{equation*}
            x_2f_{33}=x_2\psi_3+x_2^2(\partial_1\psi_3+\partial_2\psi_3)
        \end{equation*}
    \end{enumerate}
    Combining all these identities with the above expression obtained from \eqref{eq: s38} we get the identity \eqref{eq: cond for psi2}.
    Now, it's enough to insert the three expressions of $\psi_1,\psi_2,\psi_3$ in equation \eqref{eq: s32}, \eqref{eq: s33}, \eqref{eq: s35}, \eqref{eq: s36}, \eqref{eq: s37} and \eqref{eq: s38} to get the compatibility condition in \eqref{eq: condition for existence s3}, with the divisibility conditions needed to ensure that the coordinate change is well defined.
\end{proof}

\begin{lemma}
    Let $[R_{\mathfrak{s}_3},R_k]_{FN}=0$, then \eqref{eq: condition for existence s3} holds.
\end{lemma}

\begin{proof}
    The vanishing of the Fr\"olicher-Nijenhuis brackets amounts to
    \begin{align*}
    B^1_{12} &= \partial_{2}f_{11} x_{1} + (\partial_{1}f_{32} + \partial_{2}f_{12} - \partial_{2}f_{31}) x_{2} + (\partial_{1}f_{11} + \partial_{1}f_{22} - \partial_{2}f_{21}) x_{3} - f_{31} = 0 \\[5pt]
    B^1_{13} &= \partial_{3}f_{11} x_{1} + (\partial_{1}f_{11} + \partial_{1}f_{33} + \partial_{2}f_{11} + \partial_{2}f_{13} - \partial_{3}f_{31}) x_{2} + (\partial_{1}f_{23} - \partial_{3}f_{11} - \partial_{3}f_{21}) x_{3} - f_{13} - f_{21} = 0 \\[5pt]
    B^1_{23} &=\left( -\partial_{2}f_{13} + \partial_{3}f_{12}\right)x_1+ (\partial_{1}f_{12} + \partial_{2}f_{12} + \partial_{2}f_{33} - \partial_{3}f_{32}) x_{2} - (\partial_{1}f_{13} - \partial_{2}f_{23} + \partial_{3}f_{12} + \partial_{3}f_{22}) x_{3} + f_{12} - f_{22} + f_{33} = 0 \\[5pt]
    B^2_{12} &= \partial_{1}f_{22} x_{1} - (\partial_{1}f_{12} - \partial_{1}f_{32} - \partial_{2}f_{11} - \partial_{2}f_{22} + \partial_{2}f_{31}) x_{2} + \partial_{1}f_{21} x_{3} = 0 \\[5pt]
    B^2_{13} &= \partial_{1}f_{23} x_{1} - (\partial_{1}f_{13} - \partial_{1}f_{21} - \partial_{1}f_{33} - \partial_{2}f_{21} - \partial_{2}f_{23} - \partial_{3}f_{11} + \partial_{3}f_{31}) x_{2} - \partial_{3}f_{21} x_{3} - f_{21} - f_{23} = 0 \\[5pt]
    B^2_{23} &= (\partial_{1}f_{22} - \partial_{2}f_{13} + \partial_{2}f_{22} + \partial_{2}f_{33} + \partial_{3}f_{12} - \partial_{3}f_{32}) x_{2} - (\partial_{1}f_{23} + \partial_{3}f_{22}) x_{3} = 0 \\[5pt]
    B^3_{12} &= \partial_{1}f_{32} x_{1} + \partial_{2}f_{32} x_{2} + (\partial_{1}f_{31} - \partial_{1}f_{32} + \partial_{2}f_{31}) x_{3} + f_{32} = 0 \\[5pt]
    B^3_{13} &= \partial_{1}f_{33} x_{1} + (\partial_{1}f_{31} + \partial_{2}f_{31} + \partial_{2}f_{33}) x_{2} - \partial_{1}f_{33} x_{3} + f_{31} = 0 \\[5pt]
    B^3_{23} &= (\partial_{1}f_{32} + \partial_{2}f_{32}) x_{2} - (\partial_{1}f_{33} + \partial_{2}f_{33}) x_{3} + 2 \, f_{32} = 0
\end{align*}
    \begin{itemize}
        \item First of all, we prove \eqref{eq: cond s3,1}, i.e. $x_3\mid f_{32}$. To do that, consider $B^3_{12}=0$, then we can write 
        \begin{equation*}
            x_1\partial_1f_{32}+x_2\partial_2f_{32}+x_3\partial_3f_{32}-x_3\partial_3f_{32}+ {\left({\partial_{1}f_{31}} - {\partial_{1}f_{32}} + {\partial_{2}f_{31}}\right)} x_{3} + f_{32} = 0
        \end{equation*}
        Now, we use Euler's formula for homogeneous polynomials of degree $k$, and we write
        \begin{equation*}
            x_1\partial_1f_{32}+x_2\partial_2f_{32}+x_3\partial_3f_{32}=kf_{32},
        \end{equation*}
        and hence we get
        \begin{equation*}
            {\left({\partial_{1}f_{31}} - {\partial_{1}f_{32}} + {\partial_{2}f_{31}}-\partial_3f_{32}\right)} x_{3} +(k+1) f_{32} = 0
        \end{equation*}
        and in particular it follows that $x_3\mid f_{32}$.
        \item Let $f_{32}:=x_3h$, we prove \eqref{eq: cond s3,2}, i.e. $x_3\mid F$ where 
        \[
        F:=x_2 (x_1 + x_2) h + x_1 f_{22} + x_2 (f_{31} - f_{12} - f_{33}).
        \]
        Let us start by looking at $B_{12}^3=0$:
        \[
        x_1 x_3 \partial_1 h + x_2 x_3 \partial_2 h + x_3 ( \partial_1 f_{31} +  \partial_2 f_{31} - x_3 \partial_1 h) + x_3 h = 0
        \]
        \[
        \Longrightarrow \quad x_1 \partial_1 h + x_2 \partial_2 h + \partial_1 f_{31} + \partial_2 f_{31} - x_3 \partial_1 h + h = 0
        \]
        Thus, using the fact that $x_1\partial_{1}h+x_2\partial_2h\equiv_{x_3}(k-1)h$, which is Euler's formula applied to $h$, homogeneous polynomial of degree $k-1$, we get
        \begin{equation*}
        	-kh\equiv_{x_3}\partial_1 f_{31}+\partial_2f_{31}
        \end{equation*}
        Now, from $B_{13}^3=0$, the result above and Euler's formula applied to $f_{33}$:
        \[
        x_1 \partial_1 f_{33} + x_2 (\partial_1 f_{31} + \partial_2 f_{31}) + x_2 \partial_2 f_{33} - x_3 \partial_1 f_{33} = 0
        \]
        \[
        \Longrightarrow \quad k f_{33} - x_2 k h + f_{31} \equiv_{x_3} 0 
        \]
        \begin{equation}\label{eq: f31 mod x3}
        	\implies\quad f_{31} \equiv_{x_3} k (x_2 h - f_{33})
        \end{equation}
        Moreover
        \begin{equation}\label{eq: partial2f31}
        	x_1 \partial_1 f_{31} + x_2 \partial_2 f_{31} \equiv_{x_3} k f_{31} \underbrace{\implies}_{\partial_1 f_{31}\equiv_{x_3}-kh-\partial_2 f_{31}} (x_1 - x_2) \partial_2 f_{31} \equiv_{x_3} -k (f_{31} + x_1 h)
        \end{equation}
        Now, looking at $B_{12}^1$:
        \[
        x_1 \partial_2 f_{11} + x_2 (\partial_2 f_{12} - \partial_2 f_{31}) - f_{31} \equiv_{x_3} 0
        \]
        while from $B_{12}^2$:
        \[
        x_1 \partial_1 f_{22} - x_2 (\partial_1 f_{12} - \partial_2 f_{11} - \partial_2 f_{22} + \partial_2 f_{31}) \equiv_{x_3} 0
        \]
        Thus, computing $x_1 B_{12}^2 - x_2 B_{12}^1\equiv_{x_3}0$ we get
        \begin{align*}
        	 x_1^2 \partial_1 f_{22} - x_1 x_2 \partial_1 f_{12} + x_1 x_2 \partial_2 f_{11} + x_1 x_2 \partial_2 f_{22} - x_1 x_2 \partial_2 f_{31} - x_1 x_2 \partial_2 f_{11} - x_2^2 \partial_2 f_{12} + x_2^2 \partial_2 f_{31} + x_2 f_{31} \equiv_{x_3} 0
        \end{align*}
        and now, using Euler formula for $f_{22}$ and $f_{12}$ we get
        \[
        x_1 k f_{22} - x_2 k f_{12} - x_2 (x_1 - x_2) \partial_2 f_{31} + x_2 f_{31} \equiv_{x_3} 0
        \]
        Finally, using \eqref{eq: partial2f31} we get
        \begin{equation}\label{eq: star}
        	k (x_1 f_{22} - x_2 f_{12}) \equiv_{x_3} -k x_1 x_2 h - (k + 1) x_2 f_{31}.
        \end{equation}
        Now, we have
        \[
        x_2 (x_1 + x_2) h + x_1 f_{22} + x_2 (f_{31} - f_{12} - f_{33}) = x_2 (x_1 + x_2) h + (x_1 f_{22} - x_2 f_{12}) + x_2 f_{31} - x_2 f_{33}
        \]
        Multiplying by $k$ and using what we obtained above we get
        \begin{align*}
        	kF&=k  x_2 (x_1 + x_2) h + k (x_1 f_{22} - x_2 f_{12}) + k x_2 f_{31} - k x_2 f_{33}\\
        	&\underbrace{{\equiv}_{x_3}}_{\eqref{eq: star}} k x_2 (x_1 + x_2) h - k x_1 x_2 h - (k + 1) x_2 f_{31} + k x_2 f_{31} - k x_2 f_{33}\\
        	&\equiv_{x_3} k x_2^2 h - x_2 f_{31} - k x_2 f_{33}\\
        	&\underbrace{\equiv_{x_3}}_{\eqref{eq: f31 mod x3}} k x_2^2 h - k x_2 (x_2 h - f_{33}) - k x_2 f_{33} \equiv_{x_3} 0
        \end{align*}
        Thus, we proved that
        \begin{equation}
        	F\equiv_{x_3}0
        \end{equation}
        which means
        \[
        x_3 \mid  x_2 (x_1 + x_2) \frac{f_{32}}{x_3} + x_1 f_{22} + x_2 (f_{31} - f_{12} - f_{33})
        \]
        and \eqref{eq: cond s3,2} is verified.
        \item Notice that \eqref{eq: cond s3,3} is equivalent to $B^1_{12}=0$.
        \item Notice that \eqref{eq: cond s3,4} is equivalent to $x_3(B^1_{13}+B^2_{23})+x_2\,B^3_{23}=0.$
        \item Notice that \eqref{eq: cond s3,5} is equivalent to $x_3B^2_{12}+x_2\,(B^3_{12}+B^3_{23})=0.$
        \item Notice that \eqref{eq: cond s3,6} is equivalent to $
            -x_2 x_3B^1_{23}+x_3^2\, B^2_{13}+x_1 x_3\, B^2_{23}+x_2 x_3\, B^3_{13}+(x_2^2+x_1x_2)B^3_{23}=0$
        \item Notice that \eqref{eq: cond s3,7} is equivalent to $x_2B^3_{12}-x_3 B^3_{13}=0.$
        \item Notice that \eqref{eq: cond s3,8} is equivalent to $B^3_{23}=0$.
    \end{itemize}
\end{proof}

Summarizing, we have proven the following:

\begin{proposition}
    Let $\mathfrak{s}_3$ be the 3-dimensional LSA associated with linear operator
        \begin{equation*}
            R_{\mathfrak{s}_3}=\begin{pmatrix}
        x_1&-x_3&-x_2\\
        x_2&0&-x_2\\
        0&0&x_3
    \end{pmatrix}
    \end{equation*}
    Then $\mathfrak{s}_3$ is non-degenerate in the formal category.
\end{proposition}

\subsection{Non-degeneracy of $\mathfrak{s}_4$}

Let's consider the LSA $\mathfrak{s}_4$, which is associated with the Nijenhuis operator 
\begin{equation}\label{eq: op s4}
    R_{\mathfrak{s_4}}=\begin{pmatrix}
		-x_1 & x_3 & x_2 \\
		-\frac{2}{3} \, x_2 & 0 & x_3 \\
		-\frac{1}{3} \, x_3 & 0 & 0
    \end{pmatrix}
\end{equation}

\begin{lemma}
    The change of coordinates $y_i=x_i+\psi_i(x)$ with required condition \eqref{eq: condition lin} exists if and only if the components of $R_k$ satisfy the following relations:
    
\begin{subequations}\label{eq: condition for existence s4}
        \begin{empheq}[left=\empheqlbrace]{align}
        &x_3 \mid f_{32} \label{eq: cond s4,1}\\ 
        &    x_3 \mid x_1\left(x_2\dfrac{f_{32}}{x_3}-f_{22}-f_{33}\right) +\dfrac{x_2}{3}\left(f_{12}-f_{23}+\dfrac{x_2}{x_3}\left(f_{22}+2f_{33}-2\,x_2\dfrac{f_{32}}{x_3}\right)\right)\label{eq: cond s4,2}\\
        &    x_3 \mid \left(f_{22}+2f_{33}-2\,x_2\dfrac{f_{32}}{x_3}\right)\label{eq: cond s4,3} \\
            \nonumber \\
    &\begin{aligned}
    &x_{1}\left(\partial_{1}f_{22}
+\partial_{1}f_{33}\right)
+x_{2}\partial_{1}f_{31}
+\tfrac{2}{3}x_{2}
 \left(
 \partial_{2}f_{11}
 +\partial_{2}f_{22}
 +\partial_{2}f_{33}
 -\partial_{1}f_{12}
 \right)\\
    &\qquad
+\tfrac{1}{3}x_{3}
 \left(
 \partial_{3}f_{11}
 +\partial_{3}f_{22}
 +\partial_{3}f_{33}
 -\partial_{1}f_{13}
 \right)
+x_{3}\left(\partial_{1}f_{32}
+\partial_{1}f_{21}\right)
=0
    \end{aligned}
    \label{eq: cond s4,4}\\[4ex]
   &\begin{aligned}
&-\tfrac{2}{3}x_{2}\partial_{2}f_{12}
+x_{2}\partial_{2}f_{31}
-x_{3}
 \left(
 \partial_{1}f_{11}
 +\partial_{1}f_{22}
 +\partial_{1}f_{33}
 +\tfrac{1}{3}\partial_{2}f_{13}
-\partial_{2}f_{21}
-\partial_{2}f_{32}\right)\\
    &\qquad
-x_{1}\partial_{2}f_{11}
+\tfrac{1}{3}f_{12}
+f_{31}
=0
\end{aligned}
\label{eq: cond s4,5}\\[4ex]
&\begin{aligned}
&-x_{1}\partial_{3}f_{11}
-x_2\left(\tfrac{2}{3}\partial_{3}f_{12}
-\partial_{3}f_{31}
+\partial_{1}f_{11}
+\partial_{1}f_{22}
+\partial_{1}f_{33}
\right)
+x_{3}\left(\partial_{3}f_{32}
-\tfrac{1}{3}\partial_{3}f_{13}
+\partial_{3}f_{21}\right)\\
    &\qquad
-x_{3}
 \left(
 \partial_{2}f_{11}
 +\partial_{2}f_{22}
 +\partial_{2}f_{33}
 \right)
+\tfrac{2}{3}f_{13}
+f_{21}
+f_{32}
=0
\end{aligned}
\label{eq: cond s4,6}\\[4ex]
&\begin{aligned}
&2 \, {\partial_{1}f_{32}} x_{2}^{3} 
- {\left({\partial_{1}f_{13}} - 3 \, {\partial_{1}f_{21}} - 3 \, {\partial_{1}f_{32}} - {\partial_{2}f_{12}} - {\partial_{2}f_{23}} + 3 \, {\partial_{2}f_{31}}\right)} x_{3}^{3} 
+ 3 \, {\left({\partial_{1}f_{22}} + {\partial_{1}f_{33}} - {\partial_{2}f_{32}}\right)} x_{1}x_{3}^{2}\\
    &\qquad
- {\left({\partial_{1}f_{12}} - {\partial_{1}f_{23}} - 2 \, {\partial_{2}f_{11}} - {\partial_{2}f_{22}}\right)} x_{2}x_{3}^{2}
- (f_{22} - f_{33}) x_{3}^{2}
- {\left({\partial_{1}f_{22}} + 2 \, {\partial_{1}f_{33}} - 2 \, {\partial_{2}f_{32}}\right)} x_{2}^{2}x_3\\
    &\qquad
- 3 \, {\partial_{1}f_{32}} x_{1}x_2x_3
+ f_{32} x_{2}x_3=0
\end{aligned}
\label{eq: cond s4,7}\\[4ex]
&\begin{aligned}
&4 \, \partial_{2}f_{32} x_{2}^{3}
+ \left(3 \, \partial_{1}f_{11} + 3 \, \partial_{1}f_{22} + 3 \, \partial_{1}f_{33} + \partial_{3}f_{12} + \partial_{3}f_{23} - 3 \, \partial_{3}f_{31}\right) x_{3}^{3} 
+ 6 \, \partial_{1}f_{32} x_{1}x_2^2
+ 2\,f_{32} x_{2}^{2} \\
    &\qquad
+ 3 \, \left(\partial_{1}f_{12} + \partial_{1}f_{23} - 3 \, \partial_{1}f_{31} - \partial_{3}f_{32}\right) x_{1}x_3^2
+ \left(2 \, \partial_{2}f_{12} + 2 \, \partial_{2}f_{23} - 6 \, \partial_{2}f_{31} - \partial_{3}f_{22} - 2 \, \partial_{3}f_{33}\right) x_{2}x_3^2 \\
    &\qquad
- \left(f_{12} + f_{23} + 6 \, f_{31}\right) x_{3}^{2} 
- 9 \, \partial_{1}f_{32} x_{1}^{2}x_3
- 3 \, \left(\partial_{1}f_{22} + 2 \, \partial_{1}f_{33} + 2 \,\partial_{2}f_{32}\right) x_{1} x_{2}x_3\\
    &\qquad
- 2 \, \left(\partial_{2}f_{22}+ 2 \, \partial_{2}f_{33} - \partial_{3}f_{32}\right) x_{2}^{2}x_3
- 3 \, f_{32} x_{1} x_{3}=0
\end{aligned}
\label{eq: cond s4,8}\\[4ex]
&\begin{aligned}
&x_{3}^{2}\left(-3\partial_{1}f_{31}
-\partial_{2}f_{33}
+\partial_{1}f_{12}
-\partial_{2}f_{11}
-\partial_{2}f_{22}
+\partial_{1}f_{23}\right)
-x_{2}x_{3}(2\partial_{1}f_{33}
+\partial_{1}f_{22})\\
    &\qquad
-3x_{1}x_{3}\partial_{1}f_{32}
+2x_{2}^{2}\partial_{1}f_{32}
=0 \label{eq: cond s4,9}
\end{aligned}
\end{empheq}
\end{subequations}

    In this case, the coordinate change takes the form
    \begin{equation*}
        \begin{cases}\vspace{0.1cm}
            y_1=x_1 - f_{11} - f_{22} - f_{33}\\\vspace{0.1cm}
            y_2=x_2+\dfrac{1}{3} \, f_{13} - f_{21} - f_{32} + \dfrac{x_1}{x_3}\left(x_2\dfrac{f_{32}}{x_3}-f_{22}-f_{33}\right)+\dfrac{x_2}{3\,x_3}\left(f_{12}-f_{23}+\dfrac{x_2}{x_3}\left(f_{22}+2f_{33}-2\,x_2\dfrac{f_{32}}{x_3}\right)\right)\\
            y_3=x_3+\dfrac{1}{3}\left(f_{12} + f_{23}\right) - f_{31} - x_1 \dfrac{f_{32}}{x_3}+ \dfrac{1}{x_3}\left(\dfrac{x_2}{3}\left(2x_2\dfrac{f_{32}}{x_3} - f_{22} - 2f_{33}\right)\right)
        \end{cases}
    \end{equation*}
\end{lemma}

\begin{proof}
    For the operator \eqref{eq: op s4}, equation \eqref{eq: condition lin} is equivalent to
    \begin{subequations}\label{eq: conditions proof s4}
        \begin{empheq}[left=\empheqlbrace]{align}
		f_{11} + \psi_{1} &= \frac{1}{3} \, {\left(3 \, {\partial_{1}\psi_{3}} + 2 \, {\partial_{2}\psi_{1}}\right)} x_{2} + \frac{1}{3} \, {\left(3 \, {\partial_{1}\psi_{2}} + {\partial_{3}\psi_{1}}\right)} x_{3} \label{s4 eq1} \\
f_{12} - \psi_{3} &= -{\partial_{2}\psi_{1}} x_{1} + {\partial_{2}\psi_{3}} x_{2} - {\left({\partial_{1}\psi_{1}} - {\partial_{2}\psi_{2}}\right)} x_{3}  \label{s4 eq2} \\
 f_{13} - \psi_{2} &= -{\partial_{3}\psi_{1}} x_{1} - {\left({\partial_{1}\psi_{1}} - {\partial_{3}\psi_{3}}\right)} x_{2} - {\left({\partial_{2}\psi_{1}} - {\partial_{3}\psi_{2}}\right)} x_{3} \label{s4 eq3} \\
f_{21} + \frac{2}{3} \, \psi_{2} &= {\partial_{1}\psi_{2}} x_{1} - \frac{2}{3} \, {\left({\partial_{1}\psi_{1}} - {\partial_{2}\psi_{2}}\right)} x_{2} + \frac{1}{3} \, {\left(3 \, {\partial_{1}\psi_{3}} + {\partial_{3}\psi_{2}}\right)} x_{3} \label{s4 eq4} \\
f_{22} &= -\frac{2}{3} \, {\partial_{2}\psi_{1}} x_{2} - {\left({\partial_{1}\psi_{2}} - {\partial_{2}\psi_{3}}\right)} x_{3} \label{s4 eq5} \\
 f_{23} - \psi_{3} &= -\frac{1}{3} \, {\left(3 \, {\partial_{1}\psi_{2}} + 2 \, {\partial_{3}\psi_{1}}\right)} x_{2} - {\left({\partial_{2}\psi_{2}} - {\partial_{3}\psi_{3}}\right)} x_{3} \label{s4 eq6} \\
 f_{31} +\frac{1}{3} \, \psi_{3} &= {\partial_{1}\psi_{3}} x_{1} + \frac{2}{3} \, {\partial_{2}\psi_{3}} x_{2} - \frac{1}{3} \, {\left({\partial_{1}\psi_{1}} - {\partial_{3}\psi_{3}}\right)} x_{3} \label{s4 eq7} \\
 f_{32} &= -\frac{1}{3} \, {\left(3 \, {\partial_{1}\psi_{3}} + {\partial_{2}\psi_{1}}\right)} x_{3}  \label{s4 eq8} \\
f_{33} &= -{\partial_{1}\psi_{3}} x_{2} - \frac{1}{3} \, {\left(3 \, {\partial_{2}\psi_{3}} + {\partial_{3}\psi_{1}}\right)} x_{3}  \label{s4 eq9} 
        \end{empheq}
    \end{subequations}
Summing equation \eqref{s4 eq1}, \eqref{s4 eq5} and \eqref{s4 eq9} one gets the expression for $\psi_1$: 
\begin{equation*}
    \psi_1= - f_{11} - f_{22} - f_{33}.
\end{equation*}
 Computing $-\frac{1}{3}(\eqref{s4 eq2}+\eqref{s4 eq6})+\eqref{s4 eq7}$ one get: 
 \begin{equation*}
\psi_3=\frac{1}{3} (f_{12} + f_{23}) - f_{31}+x_1 \, ({\partial_{1}\psi_{3}}  + \frac{1}{3} \, {\partial_{2}\psi_{1}}) + \frac{x_2}{3}\,\left({\partial_{1}\psi_{2}} + {\partial_{2}\psi_{3}}+ \frac{2}{3} \, {\partial_{3}\psi_{1}}\right).
\end{equation*}
Using \eqref{s4 eq8} and computing  $\eqref{s4 eq5}+2\eqref{s4 eq9}-2\tfrac{x_2}{x_3}\eqref{s4 eq8}$ one get:
\begin{equation*}
    \psi_3=\frac{1}{3} (f_{12} + f_{23}) - f_{31}-\frac{f_{32}}{x_3}x_1 + \frac{x_2}{3} \,  \frac{1}{x_3}\left(2x_2\dfrac{f_{32}}{x_3} - f_{22} - 2f_{33}\right).
\end{equation*}
For $\psi_2$, one computes $\eqref{s4 eq4}-\tfrac{1}{3}\eqref{s4 eq3}$ and get
\begin{equation*}
\psi_{2} = \frac{1}{3} \, f_{13} - f_{21} + \frac{1}{3} \, {\left(3 \, {\partial_{1}\psi_{2}} + {\partial_{3}\psi_{1}}\right)} x_{1} - \frac{1}{3} \, {\left({\partial_{1}\psi_{1}} - 2 \, {\partial_{2}\psi_{2}} + {\partial_{3}\psi_{3}}\right)} x_{2} + \frac{1}{3} \, {\left(3 \, {\partial_{1}\psi_{3}} + {\partial_{2}\psi_{1}}\right)} x_{3}.
\end{equation*}
Then, using \eqref{s4 eq8}
\begin{equation*}
\psi_{2} = \frac{1}{3} \, f_{13} - f_{21} -f_{32}+ \frac{1}{3} \, {\left(3 \, {\partial_{1}\psi_{2}} + {\partial_{3}\psi_{1}}\right)} x_{1} - \frac{1}{3} \, {\left({\partial_{1}\psi_{1}} - 2 \, {\partial_{2}\psi_{2}} + {\partial_{3}\psi_{3}}\right)} x_{2}.
\end{equation*}
We compute \eqref{s4 eq5} + \eqref{s4 eq9} and get:
\begin{equation*}
    f_{22} + f_{33} = -\frac{1}{3} \, \left(3\, {\partial_{1}\psi_{3}} +2\,{\partial_{2}\psi_{1}} \right) x_{2} - \frac{1}{3} \, {\left(3 \, {\partial_{1}\psi_{2}} + {\partial_{3}\psi_{1}}\right)} x_{3}
\end{equation*}
Now we substitute the last term of the RHS in the previous expression for $\psi_2$ and using again \eqref{s4 eq8} we get
\begin{equation*}
    \psi_{2} = \frac{1}{3} \, f_{13} - f_{21} -f_{32}+ \frac{x_1}{x_3} \, {\left(x_2 \frac{f_{32}}{x_3}-f_{22}-f_{33}\right)} -\frac{x_1x_2}{3\,x_3} \partial_2\psi_1 - \frac{1}{3} \, {\left({\partial_{1}\psi_{1}} - 2 \, {\partial_{2}\psi_{2}} + {\partial_{3}\psi_{3}}\right)} x_{2}.
\end{equation*}
We compute $\eqref{s4 eq2} - \eqref{s4 eq6}$ and get:
\begin{equation*}
f_{12} - f_{23} = -{\partial_{2}\psi_{1}} x_{1} +\frac{1}{3} \, {\left(3 \, {\partial_{1}\psi_{2}} +3\,{\partial_{2}\psi_{3}}+ 2 \, {\partial_{3}\psi_{1}}\right)} x_{2}- {\left({\partial_{1}\psi_{1}} - 2\,{\partial_{2}\psi_{2}}+ {\partial_{3}\psi_{3}}\right)} x_{3} 
\end{equation*}
Now we substitute again the last term of the RHS in the last expression for $\psi_2$ we get exactly
\begin{equation*}
\dfrac{1}{3} \, f_{13} - f_{21} - f_{32} + \dfrac{x_1}{x_3}\left(x_2\dfrac{f_{32}}{x_3}-f_{22}-f_{33}\right)+\dfrac{x_2}{3\,x_3}\left(f_{12}-f_{23}+\dfrac{x_2}{x_3}\left(f_{22}+2f_{33}-2\,x_2\dfrac{f_{32}}{x_3}\right)\right)
\end{equation*}

Now, it's enough to insert the three expressions of $\psi_1,\psi_2,\psi_3$ in equation \eqref{s4 eq1}, \eqref{s4 eq2}, \eqref{s4 eq3}, \eqref{s4 eq5}, \eqref{s4 eq7} and \eqref{s4 eq8} to get the compatibility condition in \eqref{eq: condition for existence s4}, with the divisibility conditions needed to ensure that the coordinate change is well defined.
\end{proof}

Before proving that the vanishing of the Fr\"olicher-Nijenhuis brackets implies the above conditions, we recall the following technical result.

\begin{lemma}[weighted Euler Lemma]\label{lm: wEuler}
	Let $\mathcal{W}:=3x_1\partial_1+2x_2\partial_2.$ If $f\in\mathbb{R}[x_1,x_2]$ satisfies $\mathcal{W}(f)=\lambda f$ for some $\lambda \in \mathbb{R}$, then $f$ is a linear combination of the monomials $x_1^mx_2^n$ with $3m+2n=\lambda$, i.e. 
	\begin{equation*}
		f =\sum_{m,n \mid 3m+2n=\lambda}a_{m,n}x_1^mx_2^n
	\end{equation*}
\end{lemma}
\begin{proof}
	Let $f=\sum_{m,n}a_{m,n}x_1^mx_2^n$ and $\lambda \in\mathbb{R}$ such that $\mathcal{W}(f)=\lambda f$.  
	\begin{equation*}
		\mathcal{W}(f)=3x_1\partial_1f+2x_2\partial_2f=3x_1\partial_1\sum_{m,n}a_{m,n}x_1^mx_2^n+2x_2\partial_2\sum_{m,n}a_{m,n}x_1^mx_2^n = \sum_{m,n}(3m+2n)a_{m,n}x_1^mx_2^n 
	\end{equation*}
	has to be equal to $\lambda f = \sum_{m,n}\lambda a_{m,n} x_1^mx_2^n$. Then, it must hold $(3m+2n-\lambda)a_{m,n}=0$,for every monomial of $f$, which implies $3m+2n=\lambda$ as wished.\\  
\end{proof}

\begin{lemma}
    Let $[R_{\mathfrak{s}_4},R_k]_{FN}=0$, then \eqref{eq: condition for existence s4} holds.
\end{lemma}

\begin{proof}
    The vanishing of the Fr\"olicher-Nijenhuis brackets $[R_{\mathfrak{s}_4},R_k]_{FN}=:B$ amounts to
    \begin{align*}
        B^1_{12}&= -{\partial_{2}f_{11}} x_{1} - \frac{1}{3} \, {\left(3 \, {\partial_{1}f_{32}} + 2 \, {\partial_{2}f_{12}} - 3 \, {\partial_{2}f_{31}}\right)} x_{2} - \frac{1}{3} \, {\left(3 \, {\partial_{1}f_{11}} + 3 \, {\partial_{1}f_{22}} - 3 \, {\partial_{2}f_{21}} + {\partial_{3}f_{12}}\right)} x_{3} + \frac{1}{3} \, f_{12} + f_{31} = 0 \\
B^1_{13}&= -{\partial_{3}f_{11}} x_{1} - \frac{1}{3} \, {\left(3 \, {\partial_{1}f_{11}} + 3 \, {\partial_{1}f_{33}} + 2 \, {\partial_{2}f_{13}} - 3 \, {\partial_{3}f_{31}}\right)} x_{2} - \frac{1}{3} \, {\left(3 \, {\partial_{1}f_{23}} + 3 \, {\partial_{2}f_{11}} + {\partial_{3}f_{13}} - 3 \, {\partial_{3}f_{21}}\right)} x_{3} + \frac{2}{3} \, f_{13} + f_{21} = 0 \\
B^1_{23}&={\partial_{2}f_{13}} x_{1} - {\partial_{3}f_{12}} x_{1} - {\left({\partial_{1}f_{12}} + {\partial_{2}f_{33}} - {\partial_{3}f_{32}}\right)} x_{2} + {\left({\partial_{1}f_{13}} - {\partial_{2}f_{12}} - {\partial_{2}f_{23}} + {\partial_{3}f_{22}}\right)} x_{3} + f_{22} - f_{33} = 0 \\
B^2_{12}&=-{\partial_{1}f_{22}} x_{1} + \frac{2}{3} \, {\left({\partial_{1}f_{12}} - {\partial_{2}f_{11}} - {\partial_{2}f_{22}}\right)} x_{2} - \frac{1}{3} \, {\left(3 \, {\partial_{1}f_{21}} + 3 \, {\partial_{1}f_{32}} - 3 \, {\partial_{2}f_{31}} + {\partial_{3}f_{22}}\right)} x_{3} = 0 \\
B^2_{13}&=-{\partial_{1}f_{23}} x_{1} + \frac{1}{3} \, {\left(2 \, {\partial_{1}f_{13}} - 3 \, {\partial_{1}f_{21}} - 2 \, {\partial_{2}f_{23}} - 2 \, {\partial_{3}f_{11}}\right)} x_{2} - \frac{1}{3} \, {\left(3 \, {\partial_{1}f_{33}} + 3 \, {\partial_{2}f_{21}} + {\partial_{3}f_{23}} - 3 \, {\partial_{3}f_{31}}\right)} x_{3} + \frac{1}{3} \, f_{23} + f_{31} = 0 \\
B^2_{23}&=-\frac{1}{3} \, {\left(3 \, {\partial_{1}f_{22}} - 2 \, {\partial_{2}f_{13}} + 2 \, {\partial_{3}f_{12}}\right)} x_{2} + {\left({\partial_{1}f_{23}} - {\partial_{2}f_{22}} - {\partial_{2}f_{33}} + {\partial_{3}f_{32}}\right)} x_{3} + f_{32} = 0 \\
B^3_{12}&=-{\partial_{1}f_{32}} x_{1} - \frac{2}{3} \, {\partial_{2}f_{32}} x_{2} + \frac{1}{3} \, {\left({\partial_{1}f_{12}} - 3 \, {\partial_{1}f_{31}} - {\partial_{2}f_{11}} - {\partial_{3}f_{32}}\right)} x_{3} - \frac{1}{3} \, f_{32} = 0 \\
B^3_{13}&=-{\partial_{1}f_{33}} x_{1} - \frac{1}{3} \, {\left(3 \, {\partial_{1}f_{31}} + 2 \, {\partial_{2}f_{33}}\right)} x_{2} + \frac{1}{3} \, {\left({\partial_{1}f_{13}} - 3 \, {\partial_{2}f_{31}} - {\partial_{3}f_{11}} - {\partial_{3}f_{33}}\right)} x_{3} = 0 \\
B^3_{23}&=-{\partial_{1}f_{32}} x_{2} + \frac{1}{3} \, {\left(3 \, {\partial_{1}f_{33}} + {\partial_{2}f_{13}} - 3 \, {\partial_{2}f_{32}} - {\partial_{3}f_{12}}\right)} x_{3} = 0
    \end{align*}
 \begin{itemize}
        \item First of all, we prove \eqref{eq: cond s4,1}, i.e. that $x_3\mid f_{32}$. To do that, consider $B^3_{23}=0$, then we can write  
        \begin{equation*}
            {\partial_{1}f_{32}} x_{2} \equiv_{x_3} 0
        \end{equation*}
        and hence it holds ${\partial_{1}f_{32}} \equiv_{x_3}0$, which implies, looking at $B^3_{12}=0$, 
        \begin{equation*}
            {2\partial_{2}f_{32}} x_{2} +f_{32}\equiv_{x_3} 0.
        \end{equation*}
        Hence we conclude
        \begin{equation*}
            f_{32}\equiv_{x_3}0,\text{ i.e. }\quad f_{32}=x_3h
        \end{equation*}
        for a suitable $h$.
        \item Let now $f_{32}=x_3h$, we prove \eqref{eq: cond s4,3}, i.e. that $x_3\mid E$ where
        \begin{equation*}
            E:= f_{22}+2 f_{33}-2 x_2 h.
        \end{equation*}
        We start from equation $B^3_{23}$ and write: 
        \begin{equation*}
        {\partial_{3}f_{12}} -{\partial_{2}f_{13}} = 3 \left(-{\partial_1hx_2}+{\partial_{1}f_{33}}-{\partial_{2}hx_3}\right).            
        \end{equation*}
        We substitute it in $B^2_{23}$: 
        \begin{equation*}
            -\, {\partial_{1}f_{22}} x_2-2 \, {\left( -{\partial_1hx_2}+{\partial_{1}f_{33}}-{\partial_{2}hx_3}\right)} x_{2} + {\left({\partial_{1}f_{23}} - {\partial_{2}f_{22}} - {\partial_{2}f_{33}} + {\partial_{3}f_{32}}+ h\right)} x_{3} = 0 
        \end{equation*}
        and hence
        \begin{equation*}
            {\partial_{1}\left(f_{22}-2 hx_2+2f_{33}\right)} \equiv_{x_3} 0 
        \end{equation*}
        i.e.
        \begin{equation}\label{eq: s4,E0}
             \partial_1E \equiv_{x_3} 0. 
        \end{equation}
        Now we consider $B^3_{12},B^2_{12}$ and $B^3_{13}$ modulo $x_3$, which gives. 
        \begin{subequations}
        \begin{align}
        \partial_1 f_{12} - \partial_2 f_{11} - 3\partial_1 f_{31} - 2h - 3x_1 \partial_1 h - 2x_2 \partial_2 h \equiv_{x_3} 0 \label{eq: s4,E1}\\
        -x_1 \partial_1 f_{22} + \frac{2}{3}x_2 \left( \partial_1 f_{12} - \partial_2 f_{11} - \partial_2 f_{22} \right) &\equiv_{x_3} 0 \label{eq: s4,E2}\\
        3x_1 \partial_1 f_{33} + 3x_2 \partial_1 f_{31} + 2x_2 \partial_2 f_{33} &\equiv_{x_3} 0\label{eq: s4,E3}
        \end{align}
        \end{subequations}
        Computing $2x_2\,\eqref{eq: s4,E1}-3\,\eqref{eq: s4,E2}+2\,\eqref{eq: s4,E3}$ we obtain
        \begin{equation*}
            3x_1\partial_1E+2x_2\partial_2 E \equiv_{x_3}0
        \end{equation*}
        which combined with \eqref{eq: s4,E0} implies $\partial_2 E\equiv_{x_3}0$ and hence
        \begin{equation*}
            E\equiv_{x_3}0,\text{ i.e. }\quad E=x_3g
        \end{equation*}
        for a suitable $g$.
    \item Now we show \eqref{eq: cond s4,2}. Using the above notation, it can be rewritten as 
    \begin{equation*}
        x_3 \mid F, \quad\text{ where } \quad F:= x_1\left(x_2h-f_{22}-f_{33}\right) +\dfrac{x_2}{3}\left(f_{12}-f_{23}+x_2 g\right)
    \end{equation*}
    To show this condition, we use the weighted Euler Lemma \ref{lm: wEuler} and the following conditions:   
    \begin{subequations}
    \begin{align}
    B^1_{12}\equiv_{x_3}0  &\iff  -x_1 \partial_2 f_{11} - \tfrac{2}{3} x_2 \partial_2 f_{12} + x_2 \partial_2 f_{31} + \tfrac{1}{3} f_{12} + f_{31}\equiv_{x_3} 0 \label{eq: s4 F,0}\\
    B^2_{12}\equiv_{x_3}0 &\iff -x_1 \partial_1 f_{22} + \tfrac{2}{3} x_2 \left( \partial_1 f_{12} - \partial_2 f_{11} - \partial_2 f_{22} \right) \equiv_{x_3} 0 \label{eq: s4 F,1}\\
    B^2_{13}\equiv_{x_3}0 &\iff -x_1 \partial_1 f_{23} + \tfrac{1}{3} x_2 \left( 2 \partial_1 f_{13} - 3 \partial_1 f_{21} - 2 \partial_2 f_{23} - 2 \partial_3 f_{11} \right) + \tfrac{1}{3} f_{23} + f_{31}\equiv_{x_3} 0 \label{eq: s4 F,2}\\
    \partial_3 B^2_{12}\equiv_{x_3}0 &\iff -x_1 \partial_{13} f_{22} + \tfrac{2}{3} x_2 \left( \partial_{13}f_{12} - \partial_{23}f_{11} - \partial_{23}f_{22} \right) - \tfrac{1}{3} \partial_3 f_{22} - \partial_1 f_{21} + \partial_2 f_{31}\equiv_{x_3} 0 \label{eq: s4 F,3}\\
    \partial_3B^3_{13}\equiv_{x_3}0 &\iff -x_1 \partial_{13} f_{33} - \tfrac{1}{3} x_2 \left( 3 \partial_{13}f_{31} + 2 \partial_{23}f_{33} \right) + \tfrac{1}{3} \left( \partial_1 f_{13} - 3 \partial_2 f_{31} - \partial_3 f_{11} - \partial_3 f_{33} \right)\equiv_{x_3} 0 \label{eq: s4 F,4}\\
    \partial_{33}B^3_{12}\equiv_{x_3}0 &\iff -x_1 \partial_{13} h - \tfrac{2}{3} x_2 \partial_{23} h - \partial_3 h + \tfrac{1}{3} \left( \partial_{13} f_{12} - \partial_{23} f_{11} - 3 \partial_{13} f_{31} \right)\equiv_{x_3} 0 \label{eq: s4 F,5}
    \end{align}
    \end{subequations}
    It holds
    \begin{equation*}
        3x_1 \partial_1 F + 2x_2 \partial_2 F - 3F = -x_2 \,\eqref{eq: s4 F,0} + \frac{3}{2}x_1 \,\eqref{eq: s4 F,1} + x_2 \,\eqref{eq: s4 F,2} - x_2^2 \,\eqref{eq: s4 F,3} - 2x_2^2 \,\eqref{eq: s4 F,4} + 2x_2^3 \,\eqref{eq: s4 F,5} - \frac{3}{2}x_1^2 \partial_1 E - x_1 x_2 \partial_2 E
    \end{equation*}
    which means
    \begin{equation*}
        3x_1 \partial_1 F + 2x_2 \partial_2 F - 3F \equiv_{x_3} 0 
    \end{equation*}
    Now we can apply Lemma \ref{lm: wEuler} with $\lambda = 3$. This tell us that $F$ modulo $x_3$ contains only the monomial $x_1$. Moreover, $F$ is homogeneous of degree $k+1$, and hence it cannot have a part of degree 1. Hence $F\equiv_{x_3}0$. 
    \item Notice that \eqref{eq: cond s4,4} is equivalent to $-(B^3_{13}+B^2_{12})=0$.
    \item Notice that \eqref{eq: cond s4,5} is equivalent to $B^1_{12}-B^3_{23}=0$.
    \item Notice that \eqref{eq: cond s4,6} is equivalent to $B^1_{13}+B^2_{23}=0$.
    \item Notice that \eqref{eq: cond s4,7} is equivalent to $-x_3^2 B^1_{23}- 3x_3^2 B^2_{12}+x_2x_3 B^2_{23}+(3x_1x_3-2x_2^2)B^3_{23}=0$.
    \item Notice that \eqref{eq: cond s4,8} is equivalent to $-3x_3^2(B^1_{12}+B^2_{13})+3x_2x_3(2B^3_{13}+B^2_{12})+9x_1x_3B^3_{12}-6x_2^2B^3_{12}=0$.
    \item Notice that \eqref{eq: cond s4,9} is equivalent to $x_3B^2_{23}-2x_2B^3_{23}+3x_3B^3_{12}=0$.
\end{itemize}
\end{proof}

Summarizing, we have proven the following:

\begin{proposition}
    Let $\mathfrak{s}_4$ be the 3-dimensional LSA associated with linear operator
        \begin{equation*}
        R_{\mathfrak{s}_4}=\begin{pmatrix}
		      -x_1 & x_3 & x_2 \\
		      -\frac{2}{3} \, x_2 & 0 & x_3 \\
		      -\frac{1}{3} \, x_3 & 0 & 0
        \end{pmatrix}
        \end{equation*}
    Then $\mathfrak{s}_4$ is non-degenerate in the formal category.
\end{proposition}

\end{document}